\documentclass[12pt,twoside]{article}
\usepackage{amsmath}

\newcommand{\beq}{\begin{equation}}
\newcommand{\enq}{\end{equation}}

\newcommand{\Rb}{\mathbb{R}}
\newcommand{\Rbn}{\mathbb{R}^{n-1}}

\newcommand{\pip}{\Pi_{+}}

\newcommand{\om}{\omega}
\newcommand{\f}{\frac}
\newcommand{\pa}{\partial}
\newcommand{\qed}{\mbox{$\quad\blacksquare$}}

\newcommand{\lam}{\lambda}

\newcommand{\la}{\lambda}

\newcommand{\cl}{C^{\lambda}}

\newcommand{\Ff}{F_{\lam, M}[f]}

\newcommand{\IN}{\infty}

\usepackage{amssymb}

\def\ges{\,\,\lower1ex\hbox{$\stackrel{\tss >}{\sim}$}\,\,}
\def\les{\,\,\lower1ex\hbox{$\stackrel{\tss <}{\sim}$}\,\,}

\newcommand{\be}{\begin{equation}}
\newcommand{\ee}{\end{equation}}
\newcommand{\dss}{\displaystyle}
\newcommand{\tss}{\textstyle}
\newcommand{\veps}{\varepsilon}

\newcommand{\ol}{\overline}

\newcommand{\raro}{\rightarrow}

\newtheorem{lem}{Lemma}[section]
\newtheorem{cor}{Corollary}[section]
\newtheorem{theo}{Theorem}[section]
\newtheorem{defin}{Definition}[section]
\newtheorem{prop}{Proposition}[section]
\newtheorem{rem}{Remark}[section]

\numberwithin{equation}{section}

\begin{document}

\raisebox{7ex}[1ex]{\fbox{{\footnotesize To appear in 
{\it Potential analysis}}}}

\begin{center}
{\large {\bf Sharp Growth Estimates for Modified Poisson
Integrals\\
in a Half Space}}\\
\medskip

David Siegel\footnote{
Research partially supported by an NSERC Individual Research Grant.}
 (dsiegel@math.uwaterloo.ca)\\
Department of Applied Mathematics\\
University of Waterloo\\
Waterloo, Ontario, Canada N2L 3G1\\
\medskip
Erik Talvila\footnote{Research partially supported by
an Ontario Graduate Scholarship.
\newline
1991 {\it Mathematics Subject Classification.}  Primary
31B10, Secondary 35J05.\newline
{\it Key words and phrases.}  Poisson integral, half space Dirichlet problem,
half space Neumann problem.}
(etalvila@math.alberta.ca)\\
Department of Mathematical Sciences\\
University of Alberta\\
Edmonton, Alberta, Canada T6G 2E2\\
\medskip
November 22, 2000
\end{center}
\bigskip

\begin{quotation}\noindent
{\bf {A\sc{bstract}.}} For continuous boundary data, including data
of polynomial growth, modified
Poisson integrals are used to write solutions to the half
space Dirichlet and Neumann problems in $\mathbb{R}^{n}$.
Pointwise growth estimates for these integrals are given and the estimates
are proved sharp in a strong sense.  For decaying data, a new type of modified
Poisson integral is introduced and used to develop asymptotic 
expansions for solutions of  these half space problems. 
\end{quotation}
\bigskip

\section{\!\!\!\!\!\!.\; Introduction.}

For $x$ in the half space $\Pi_{+} = \{ x \in \mathbb{R}^{n}
| x_{n} > 0 \}\,\, (n \geq 2)$ let $y \in \mathbb{R}^{n-1}$ be identified
with the projection of $x$ onto the hyperplane $\partial
\Pi_{+}$.  Let $B_{r}(a)$ be the open ball in
$\mathbb{R}^{n-1}$ with centre $a \in \mathbb{R}^{n-1}$ and radius $r
> 0$.  When $a=0$ we write $B_{r}$.  And let $\theta$ be
the angle between $x$ and $\hat{e}_{n}$, i.e., $x_{n}
= |x| \cos \theta$, $| y |= |x
|\sin \theta$ and $0 \leq \theta < \pi /2$ when $x
\in \Pi_{+}$.  We will write $x = {\dss\sum_{i=1}^{n-1}}
y_{i} \hat{e}_{i} + x_{n} \hat{e}_{n}$ where
$\hat{e}_{i}$ is the $i^{\,\rm th}$ unit coordinate
vector and $\hat{e}_{n}$ is normal to $\partial
\Pi_{+}$.  Unit vectors will be denoted with a caret,
$\hat{x} = x/|x |$ for $x \neq 0$.

For $\lambda > 0$ $( \lambda \in \mathbb{R} )$ and $y' \in
\mathbb{R}^{n-1}$ define the kernel
\begin{equation}
K (\lambda, x, y') = \Bigl[ |y' - y |^{2} +
x_{n}^{2} \Bigr]^{-\lambda} . \label{1.1}
\end{equation}
The Poisson integrals for the half space problem $\Delta
u = 0$ $(x \in \Pi_{+})$ with Dirichlet and Neumann data
$f\!: \mathbb{R}^{n-1} \raro \mathbb{R}$ on $\partial \Pi_{+}$ are,
respectively,
\begin{equation}
D [f] (x) = \alpha_{n} x_{n} \int\limits_{\mathbb{R}^{n-1}} f(y')
K \left( \tfrac{n}{2}, x, y' \right)\, dy' \quad (n
\geq 2) \label{1.2}
\end{equation}
and
\begin{equation}
N[f] (x) = \frac{\alpha_{n}}{n-2} \int\limits_{\mathbb{R}^{n-1}}
f(y') K \left( \tfrac{n-2}{2}, x,y' \right)\, dy'
\quad (n \geq 3) . \label{1.3}
\end{equation}
Here $\alpha_{n} = 2/(n\omega_{n})$ and $\omega_{n} =
\pi^{n/2}/\Gamma (1+n/2)$ is the volume of the unit $n$-ball
.  When $n=2$, $N$ has a logarithmic kernel and is
not dealt with here.  (See the remarks at the end of
this section.)

The functions defined by \eqref{1.2} and \eqref{1.3}
will be harmonic in $\Pi_{+}$ if
\begin{equation}
\int\limits_{\mathbb{R}^{n-1}} \frac{| f(y') |\, dy'}{| y'
|^{2\lambda} +1} < \infty \label{1.4}
\end{equation}
with $\lambda = n/2$ and $(n-2)/2$, respectively
(\cite{Flett}, Theorem 6).  If $f$ is continuous then
convergence of the appropriate integral in \eqref{1.4}
is sufficient for \eqref{1.2} or \eqref{1.3} to be a
classical solution of the respective Dirichlet or
Neumann problem on $\Pi_{+}$ (cf.\ Corollary
\ref{cor2.1} and \ref{cor2.2} below).  Notice that since
$\Pi_{+}$ is unbounded the integral of $f$ over
$\mathbb{R}^{n-1}$ need not vanish for $N[f]$ to be a solution of
the Neumann problem.

When the integral in \eqref{1.4} diverges but
\begin{equation}
\int\limits_{\mathbb{R}^{n-1}} \frac{| f(y') |\, dy'}{| y'
|^{M+2 \lambda} +1} < \infty \label{1.5}
\end{equation}
for a positive integer $M$ we can use the modified kernel
\begin{equation}
K_{M} (\lambda , x, y') = K(\lambda , x, y') -
\sum_{m=0}^{M-1} \frac{| x |^{m}}{| y'
|^{m+2 \lambda}} C_{m}^{\lambda} (\sin \theta \cos
\theta' ) \label{1.6}
\end{equation}
(defined for $| y' | > 0$) where $0 \leq \theta'
\leq \pi$ is the angle between $y$ and $y'$, i.e., ${y
\!\cdot\! y'} = | y' | \, | x | \sin \theta \cos
\theta'$ and $K_{0} = K$.  If $y = 0$ or $y' = 0$ we
take $\theta' = \pi /2$.  When $n=2$, we take $\theta' =
0$ or $\pi$ according as $y'$ and $x_{1}$ are on the
same or opposite side of the origin.  Equivalently,
$\cos \theta' = {\rm sgn} (x_{1}\, y')$.  In \eqref{1.6} the
first $M$ terms of the asymptotic expansion of $K$
in inverse powers of $| y' |$ are removed.  The
coefficients are in terms of Gegenbauer polynomials,
$C_{m}^{\lambda}$, most of whose properties used herein
are derived in \cite{Szego}.

Let $w \!: \mathbb{R}^{n-1} \raro [0,1]$ be continuous so that
$w(y)\equiv 0$ for $0 \leq | y | \leq 1$ and $w(y)
\equiv 1$ for $| y | \geq 2$.  Define modified
Dirichlet and Neumann integrals
\begin{equation}
D_{M} [f] (x) = \alpha_{n} x_{n} \int\limits_{\mathbb{R}^{n-1}}
f(y') K_{M} \left( \tfrac{n}{2}, x, y' \right)\, dy' \quad
(n \geq 2) \label{1.7}
\end{equation}
\begin{equation}
N_{M} [f] (x) = \frac{\alpha_{n}}{n-2}
\int\limits_{\mathbb{R}^{n-1}} f (y') K_{M} \left(
\tfrac{n-2}{2}, x, y' \right)\, dy' \quad (n \geq 3).
\label{1.8}
\end{equation}
Then $u(x) = D_{M} [wf](x) + D[(1-w)f](x)$ and $v(x) =
N_{M} [wf](x) + N[(1-w)f](x)$ are respective solutions
of the classical half space Dirichlet and Neumann
problems.  The purpose of the function $w$ is merely to avoid
the singularity of the modified kernel at the origin.  The 
Dirichlet version of $K_{M}$ appears in
\cite{Armitage2}, \cite{SiegelTalvila} and \cite{Yoshida}, with
inspiration from \cite{FinkelsteinScheinberg}.  The
Neumann version is discussed by Gardiner
(\cite{Gardiner}) and Armitage (\cite{Armitage}).

In this paper we give growth estimates for $u$ and $v$
under \eqref{1.5} and prove they are sharp.  This is
done in Theorem \ref{theo2.1} by first defining
\begin{equation}
F_{\lambda ,M}[f] (x) = \int\limits_{| y' | > 1}
f(y') K_{M} (\lambda , x, y')\, dy' \label{1.9}
\end{equation}
and proving that
\begin{equation}
F_{\lambda , M}[f] (x) = o (| x |^{M} \sec^{2
\lambda} \theta ) \quad \text{as} \quad | x| \raro
\infty \quad \text{with} \quad x \in \Pi_{+} .  \label{1.10}
\end{equation}
The order relation is interpreted as $\mu (r)r^{-M} \raro
0$ as $r \raro \infty$ where $\mu (r)$ is the supremum
of $| F_{\lambda , M} [f](x) | \cos^{2 \lambda}
\theta$ over $x \in \Pi_{+}$, $| x | = r$.
A growth condition $\omega$ is said to be sharp if given any
function $\psi = o(\omega)$ and any sequence $\{x_{i}\} \in
\Pi_{+}$ with $|x_{i}| \to \infty$, we can find data $f$ so
that the solution corresponding to $f$ is not  $o(\psi)$  
on this sequence (see Definition \ref{defin2.1} below).  The sharpness proof
is complicated by the fact that the modified kernels are not of
one sign.  For each $x_{i}$, regions in $\mathbb{R}^{n-1}$ are   
determined where the kernel is of one sign.  Data is then
chosen so that the contribution from integrating where the sign
of the kernel is not known is cancelled out and the main contribution
comes from integrating over a neighbourhood of the singularity of the 
kernel. This proof makes up a substantial portion of the paper.
Note that condition \eqref{1.5}
is necessary and sufficient for $F_{\lambda,M}[f](x)$ to
exist  as a Lebesgue integral on $\pip$.  See Proposition 3.4.1
in \cite{talvila}.

In the third section of the paper, the modified Neumann operator
$N_{M}$ is represented as various integrals of the
modified Dirichlet operator $D_{M}$.

In the final section a second type of modified kernel is
defined, $\eqref{4.1}$, useful when
\begin{equation}
\int\limits_{\mathbb{R}^{n-1}} | f(y') |\, (| y'
|^{M-1} +1)\, dy' < \infty \label{1.11}
\end{equation}
for a positive integer $M$.  This new kernel will be used
in Theorem \ref{theo4.2} to derive asymptotic expansions
of $D[f]$ and $N[f]$ under \eqref{1.11}.  A growth estimate is
given for the remainder term and the estimate is proved sharp
as in the previous case.

The half plane $(n=2)$ Neumann Poisson integral has a
logarithmic kernel.  A modified kernel can be defined as
in \eqref{1.6} and is in some sense the limit as
$\lambda \raro 0^{+}$.  However, this case is
sufficiently different to warrant separate exposition.
Analogues of the results in this paper will be presented
elsewhere.

\section{\!\!\!\!\!\!.\; First type of modified kernel.}

The expansion \eqref{1.6} arises from the generating
function for Gegenbauer polynomials (\cite{Szego},
4.7.23)
\begin{equation}
(1-2tz+z^{2})^{-\lambda} = \sum_{m=0}^{\infty} z^{m}
C_{m}^{\lambda} (t), \quad \lambda > 0 , \label{2.1}
\end{equation}
where $C_{m}^{\lambda} (t) =
\dfrac{1}{m!}\,\dfrac{\partial^{m}}{\partial z^{m}}
(1-2tz+z^{2})^{-\lambda} \biggl|_{z=0}$. If $-1 \leq t
\leq 1$ the series converges absolutely for $| z |
< 1$ (the left side of \eqref{2.1} is singular at $z = t
\pm i \sqrt{1-t^{2}}$).  The majorisation and derivative
formulas
\begin{align}
| C_{m}^{\lambda} (t) | & \leq C_{m}^{\lambda} (1)
= \binom{2\lambda + m -1}{m} = \frac{\Gamma(2\lambda
 +m)}{\Gamma(2\lambda)\,\Gamma(m+1)}
\label{2.2} \\
\dfrac{d}{dt} C_{m}^{\lambda} (t) & = 2 \lambda\,
C_{m-1}^{\lambda +1} (t) \label{2.3}
\end{align}
are proved in \cite{Szego} (4.7.3, 7.33.1, 4.7.14).
Hence, the series in \eqref{2.1} converges if $| z
| < 1$, uniformly for $-1 \leq t \leq 1$ and the same
can be said for all of its derivatives with respect to
$z$ and $t$.  From the definition above and Fa\`{a} di
Bruno's formula for the $m^{\rm th}$ derivative of a
composite function (\cite{AbramowitzStegun}, p. 823) it
can be seen that $C_{m}^{\lambda}(t)$ is a polynomial
in $t$ of degree $m$.  And,
\begin{equation}
C_{0}^{\lambda} (t) = 1, \quad C_{1}^{\lambda} (t) = 2
\lambda t . \label{2.4}
\end{equation}

A proof of the following lemma is hinted at in
\cite{Armitage2}, \cite{Gardiner} and
\cite{Yoshida} by
reference to a more general result on axial polynomials
in \cite{Kuran} (Theorem 2).  However, we give a simple
direct proof.

\begin{lem}\label{lem2.1}
For $m = 0,1,2,3, \ldots$ the functions
$h_{m+1}^{(0)} (x) = x_{n} | x |^{m}$
$C_{m}^{n/2}(\Theta ) (n \geq 2)$ and $h_{m}^{(1)} (x)
= | x |^{m} C_{m}^{(n-2)/2} (\Theta )$ $(n \geq
3)$
are homogeneous harmonic polynomials of
degree $m+1$ and $m$, respectively, where $\Theta =
\sin \theta \cos \theta'$.
\end{lem}

\bigskip
\noindent
{\it Proof}: Using \eqref{2.1} we obtain the expansion
of the fundamental solution of Laplace's equation
\begin{equation}
| x - x' |^{2-n} = \sum_{m=0}^{\infty} \frac{|
x |^{m}}{| x' |^{m+n-2}} C_{m}^{(n-2)/2}
(\hat{x} \!\cdot\! \hat{x}'), \quad n \geq 3 . \label{2.5}
\end{equation}
If $x' \neq 0$ this series converges for $| x | <
| x' |$ and defines a harmonic function.  Each
term is homogeneous in $x$ of degree $m$ and it is
clear from \eqref{2.4}  that the first two terms are
harmonic.  Given $x$, take $x'$ such that $| x' | >
| x |$.  Differentiating termwise in $x$ gives
\begin{equation}
\Delta | x - x' |^{2-n} = 0 = \sum_{m=2}^{\infty}
| x' |^{-(m+n-2)} \Delta \left( | x |^{m}
C_{m}^{(n-2)/2} (\hat{x} \!\cdot\! \hat{x}') \right) .
\label{2.6}
\end{equation}
Each term $\Delta \left( | x |^{m} C_{m}^{(n-2)/2}
(\hat{x} \!\cdot\! \hat{x}') \right)$ is homogeneous of
degree $m-2$, hence, by the linear independence of
homogeneous functions, $| x |^{m} C_{m}^{(n-2)/2}
(\hat{x} \!\cdot\! \hat{x}')$ is harmonic on $\mathbb{R}^{n}$ for
each $m \geq 0$.  Every harmonic function can be
uniquely written as a sum of homogeneous harmonic
polynomials so $| x |^{m}\, C_{m}^{(n-2)/2} (\hat{x}
\!\cdot\! \hat{x}')$ is a homogeneous harmonic polynomial of
degree $m$\,(\cite{Axler}, 1.26, 1.27).

Now set $x_{n}' = 0$, then $\hat{x} \!\cdot\! \hat{x}' =
y \!\cdot\! y'/(| x | \, | y' |) = \sin
\theta\,\, \hat{y} \!\cdot\! \hat{y}' = \Theta$.  Hence, the
Neumann half space expansion is
\begin{equation}
\left[ | y ' - y |^{2} + x_{n}^{2}
\right]^{-(n-2)/2} = \sum_{m = 0}^{\infty} \frac{| x
|^{m}}{| y' |^{m+n-2}} C_{m}^{(n-2)/2} (\Theta
), \quad n \geq 3 , \label{2.7}
\end{equation}
and each term in the series is a homogeneous harmonic polynomial of
degree $m$.

For the Dirichlet expansion differentiate \eqref{2.5}
with respect to $x_{n}'$, use \eqref{2.3} and
\eqref{2.4}, and set $x_{n}' = 0$.   Then for $n \geq 3$
\begin{equation}
x_{n} \left[ | y' - y |^{2} + x_{n}^{2}
\right]^{-n/2} = \sum_{m=0}^{\infty} \frac{x_{n}| x
|^{m}}{| y' |^{m+n}} C_{m}^{n/2} (\Theta ),
\label{2.8}
\end{equation}
and each term in the series is a homogeneous harmonic
polynomial of degree $m+1$.

When $n =2$, use $C_{m}^{1} (\cos \phi
) = \sin [(m+1) \phi ]  \csc\phi$.  Then from
\eqref{2.5} we recover the trigonometric expansion
\begin{equation}
\frac{x_{2}}{(\xi -x_{1})^{2} +x_{2}^{2}} =
\sum_{m=1}^{\infty} \frac{r^{m} \sin (m
\phi)}{\xi^{m+1}}, \quad r < |\xi|, \label{2.9}
\end{equation}
where we have written $r= | x |$, and $\phi =
\pi/2 - \theta$ to conform with the usual polar
coordinates $(x_{1} = r \cos \phi, \, x_{2} = r \sin \phi)$.
Each $r^{m} \sin (m\phi)$ is a homogeneous harmonic
polynomial of degree $m$.  \qed

\begin{rem}
\label{rem2.2}{\rm When $x_{n} = 0$, $h_{m}^{(0)}$ and
$\partial h_{m}^{(1)}/\partial x_{n}$ vanish.  The
spherical harmonics of degree $m$ are the restriction of
the homogeneous harmonic polynomials to the unit
sphere.  If we write $Y_{m}^{(0)} (\hat{x}) =
h_{m}^{(0)} (\hat{x})$ and $Y_{m}^{(1)} (\hat{x}) =
h_{m}^{(1)} (\hat{x})$ then $h_{m}^{(0)}(x) = | x
|^{m} Y_{m}^{(0)} (\hat{x})$ and $h_{m}^{(1)} (x) =
| x |^{m} Y_{m}^{(1)} (\hat{x})$.  The functions
$| x |^{-(m+n-2)}Y_{m}^{(0)} (\hat{x})$ and $|
x |^{-(m+n-2)} Y_{m}^{(1)}(\hat{x})$ are harmonic for
$| x | > 0$ (interchange $x$ and $x'$ in
\eqref{2.5} and \eqref{2.6}).} 
\end{rem}

\bigskip
Our first theorem will be a sharp growth estimate for
$F_{\lambda , M}$.  First we introduce the following
definition.

\begin{defin}
\label{defin2.1}
Let $\omega \!: \Pi_{+} \raro (0, \infty )$ then $\omega$
is a {\it sharp} growth condition for $F_{\lambda , M}$
if
\begin{itemize}
\item[(i)]
$F_{\lambda , M} [f](x) = o (\omega (x))$ $(x \in
\Pi_{+}, \, | x | \raro \infty )$, for all $f$
satisfying \eqref{1.5}
 \item[(ii)]
If $\psi\! \!: \Pi_{+} \raro (0,\infty )$ and $\psi (x) = o
(\omega (x))$ then for any sequence $\{ x^{(i)}\}$ in
$\Pi_{+}$ such that $| x^{(i)} | \raro \infty$ as
$i \raro \infty$ there exists a continuous function $f$
satisfying \eqref{1.5} with $
F_{\lambda , M} [f] (x^{(i)})/ \psi (x^{(i)}) \not \raro 0$
as $i \raro \infty$.
\end{itemize}
\end{defin}
Note that it is essential that the limit condition on $F_{\la,M}[f]/\psi$
be checked on all paths to infinity.  For example, $\omega_1(x)=|x|$
and $\omega_2(x)=|x|\sec\theta$ agree on all radial paths but allow
very different behaviour on paths approaching $\pa\pip$.

\bigskip
Let
\begin{equation}
\Phi_{\pm} (\Theta , \zeta ) = M C_{M}^{\lambda} (\Theta
) \pm (2 \lambda + M-1) C_{M-1}^{\lambda} (\Theta )\, \zeta
\quad \text{and} \quad \Theta = \sin \theta \cos \theta'.
\label{2.10}
\end{equation}
The integral representation for $M \geq 1$
\begin{equation}
K_{M} (\lambda , x, y') = K(\lambda , x, y')
\int\limits_{\zeta =0}^{| x | / | y' |} (1-2
\Theta \zeta + \zeta^{2})^{\lambda-1} \Phi_{-}(\Theta ,
\zeta )\, \zeta^{M-1}\, d \zeta \label{2.11}
\end{equation}
was derived in \cite{SiegelTalvila} (Theorem 5.1) by
summing a Gegenbauer recurrence relation.  Use of
\eqref{2.11} allows us to prove

\begin{theo}
\label{theo2.1}
Let $\lambda > 0$ and $f$ be measurable so that
\eqref{1.5} holds for integer $M \geq 0$.  Then
$F_{\lambda , M}[f](x) = o \left( | x |^{M}
\sec^{2\lambda} \theta \right)$ $(x \in \Pi^{+}, | x
|
\raro \infty )$ and the order relation is sharp in the
above sense.
\end{theo}

\bigskip
A weaker form of sharpness (with respect to the
exponents of $| x |$ and $\sec \theta$) was
obtained for $D[f]$ in \cite{SiegelTalvila}.  Also, a
growth estimate like \eqref{1.10} was obtained for $D_{M}
[f]$ but we provide a shorter proof here.  Sharpness of
the order relation is proven by finding regions in
$\mathbb{R}^{n-1}$ where $K_{M}(\lambda, x, y')$ is of
one sign (for fixed $x$).  Data $f$ is then chosen large
enough so that $F_{\lambda,M}(x)$ is  positive for all
$x$ and equal to $\psi(x)$ on a subsequence of the given
sequence in (ii) of Definition \ref{defin2.1}.  The proof is quite long but has been broken down into digestible
pieces as detailed below.  Of crucial importance is the integral
form of the modified kernel, given in \eqref{2.11}.

\medskip
\noindent
{\bf Step I}\,\,
It is shown that $\Ff =o(|x|^M\sec^{2\la}\theta)$ for any
measurable function $f$ satisfying the integrability condition 
\eqref{1.5}.  In \eqref{2.11}, the original kernel $K$ is
estimated as was done in Theorem~2.1  of \cite{SiegelTalvila} ($\alpha =n/2$).  The
function $\Phi(\Theta,\zeta)$  has a simple zero precisely where
$1-2\theta\zeta +\zeta^2$ vanishes, at $\Theta=\zeta=1$.  So the ratio
$\Phi(\Theta,\zeta)/\sqrt{1-2\theta\zeta +\zeta^2}$ is bounded and
the integrand in \eqref{2.11} is continuous for $\la \geq 1/2$ and
unbounded but integrable when $0<\la <1/2$.  In either case, elementary
approximations lead to an upper bound for $|K_M|$ on which the 
Dominated Convergence Theorem can be used to prove \eqref{1.10}. 

\smallskip
\noindent
{\bf Step II}\,\,  
The order estimate $o(|x|^M\sec^{2\la}\theta)$ is
now proven to be sharp, first for  given sequences which have a 
subsequence $\tilde{x}^{(i)} = a_{i} \hat{e}_{1} + b_{i}
\hat{e}_{n}$ that stays bounded away from the $\hat{e}_{1}$ axis of
$\pa\pip$ by an angle $\theta_0$ ($0\leq \theta_0 <\pi/2$).  On such
a sequence the growth condition reduces to $o(|x|^M)$.  A region
$\Omega_{1} \subset \mathbb{R}^{n-1}$ is found on which $\Phi_-$ and
hence $K_M$ are of one sign.  Due to the parity of $\cl_m$ about
zero ($\cl_m$  is even if $m$ is even and odd if $m$ is odd) it 
follows that $\Phi_-(\Theta,\zeta)$ will be of one sign if $|\Theta|$
is small enough.   Since $\Theta=\sin\theta\,\cos\theta'$, this
is accomplished by restricting $\theta'$ to lie near $\pi/2$.  And, 
$\Omega_{1}$ is taken as the region between two cones, both of which
have an opening angle of nearly $\pi/2$ from the $\hat{e}_{1}$ axis.
For $y' \in \Omega_{1}$, the combination $(-1)^{\lceil M/2\rceil}
\Phi_-(\Theta,\zeta)$ is strictly positive when $\zeta >0$.  (When
$x\in \Rb$, the {\it ceiling} of $x$, $\lceil x\rceil$, is $x$ if $x\in
\mathbb Z$ and is the next largest integer if $x\not\in\mathbb Z$.)
A lower bound on $(-1)^{\lceil M/2\rceil}K_M$ is now obtained, equation
\eqref{4.2.24}.   Data is then chosen that has support in $\Omega_{1}$
and is large on a sequence of unit half balls along the $\hat{e}_{2}$ axis.
(This is an axis orthogonal to $\hat{e}_{1}$.  Something slightly
different is done when $n=2$.)  By taking  the data large enough
we have $\limsup F_{\lambda ,M}[f]/\psi
 \geq 1$ on a subsequence and 
sharpness of the growth estimate for
this special type of sequence now follows.

\smallskip
\noindent
{\bf Step III}\,\,  Now considered are sequences with a subsequence
$\tilde{x}^{(i)} = a_{i} \hat{e}_{1} + b_{i}
\hat{e}_{n}$ that approaches the boundary at the $\hat{e}_{1}$ axis.
Again, a region is found where $\Phi_-$ is of one sign.  On the
sequence, we have $\sin\theta\to 1$ so taking $\theta'$ near $0$ 
makes $\Theta$ nearly equal to $1$.  In this case then, the kernel
$K(\la,x,y')$ will be singular for $|y'|=|x|$ and $\Theta\to 1$.
Hence, in \eqref{1.6} it will dominate the Gegenbauer terms 
subtracted from it.  A region $\Omega_{2} \subset \mathbb{R}^{n-1}$ is
defined to be the portion of a cone with $|y'|>1$ and axis along
$\hat{e}_{1}$.  
The opening angle $\theta'$ is taken small enough
so that when $y'_1 >0$, $|y'|$ is near $|x|$ and $|x|/A <|y'| <A|x|$ for
a constant $A>1$, we have $K_M>0$, i.e., near the
singularity of $K$.  The modified kernel is also positive for large
values of $|y'|$ in $\Omega_{2}$ but changes sign when $y'_1 >0$ and 
$1<|y'|<|x|/A$.
And, due to the parity of
$\cl_m$, the modified kernel is one sign when $y'\in\Omega_{2}$
with $y'_1 <0$.  Data is chosen to have support within $\Omega_{2}$
on a sequence of balls along the $\hat{e}_{1}$ axis.  When $y'_1$ is
positive, $f(y')$ is positive and when $y'_1$ is negative, $fK_M$ is
positive.  Contributions to $\int_{| y' | > 1}
f(y') K_{M} (\lambda , x, y')\, dy'$ are now known to be
positive except when integrating over $\Omega_{>}$, that portion
of $\Omega_2$ with $1<|y'|<|x|/A$.  But $f$ is
chosen so that if the reflection of $y'$ across the $y'_1=0$
hyperplane is denoted $y^*$, then if $y'_1>0$ we have 
$f(y^*)=(-1)^M A_\la f(y')$, where $A_\la >1$  is a constant.  The
data is said to be given a ``super odd" or ``super even" extension from $y'_1 >0$
to $y'_1 <0$, according as $M$ is even or odd.  This allows the
contribution from integrating over $\Omega_{>}$, where $fK_M$ is
not of one sign, to be balanced out by the contribution from
integrating over the reflection of $\Omega_{>}$ to $y'_1<0$, where
$fK_M$ is positive.  The contribution to $\int_{| y' | > 1}
f(y') K_{M} (\lambda , x, y')\, dy'$ from integrating near the
singularity of $K_M$, i.e., over $\Omega_{2}$, produces a lower
bound for $F_{\la,M}[f]$ from which it follows that $F_{\la,M}[f](x^{(i)})
/\psi(x^{(i)}) \not\to 0$, where $\psi$ and $x^{(i)}$ are given
in the theorem.  Note that all the $\Omega$ regions defined here
depend on $|x|$.

\smallskip
\noindent
{\bf Step IV}\,\,  The special case of sequences 
$\tilde{x}^{(i)} = a_{i} \hat{e}_{1} + b_{i}
\hat{e}_{n}$ considered in II and III is shown to be applicable to
general sequences in $\pip$.  Since $\pa B_+$ is compact, for any
sequence $r_i{\hat r}_i$ in $\pip$, the sequence $\{{\hat r}_i\}$
has a limit point ${\hat r}_0 \in \pa {\overline B}_+$.  This direction
is then rotated to correspond to ${\hat y}_1$.

\bigskip
\noindent
{\it Proof}:   Write $s = | x | / | y' |$.
Throughout the proof $d_{1}, d_{2}, \ldots, d_{9}$ will be
positive constants (depending on $\lambda$ and $M$).

\smallskip
\noindent
{\bf Step I}\,\,  First suppose $M \geq 1$.

In \cite{Szego} (4.7.27) for $M \geq 2$ we have
\begin{equation}
MC_{M}^{\lambda} (t) = (2 \lambda + M-1)\, t\,
C_{M-1}^{\lambda} (t) - 2 \lambda\, (1-t^{2})\,
C_{M-2}^{\lambda +1} (t) . \label{4.2.12}
\end{equation}
With reference to \eqref{2.10} and \eqref{2.2} we can write
\begin{align}
\frac{| \Phi_{-}(\Theta , \zeta )|}{\sqrt{1-2
\Theta \zeta + \zeta^{2}}} &= \frac{\left| (2 \lambda +
M-1) (\Theta - \zeta ) C_{M-1}^{\lambda} (\Theta ) - 2
\lambda (1-\Theta^{2}) C_{M-2}^{\lambda +1} (\Theta )
\right|}{\sqrt{(\Theta - \zeta)^{2} + (1-\Theta^{2})}}
\label{4.2.12.5} \\
& \leq (2 \lambda + M-1) \binom{ 2 \lambda + M-2}
{M-1} + 2 \lambda \binom{ 2 \lambda
+ M-1}{ M-2} \notag \\
& = 2 \lambda \binom{2 \lambda + M}{M-1}
\label{4.2.13}
\end{align}
for $M \geq 2$.  If we define $C_{-m}^{\lambda} = 0$ for
$m = 1,2, 3, \ldots$ and use the fact that
$C_{0}^{\lambda} (\Theta ) =1$ and $C_{1}^{\lambda}
(\Theta ) = 2 \lambda \Theta$ then \eqref{4.2.12} and
\eqref{4.2.13} still hold when $M=1$.  Hence, \eqref{2.11}
and \eqref{4.2.13} give
\begin{equation}
\left| K_{M} (\lambda , x, y') \right| \leq d_{1}\,
K(\lambda , x, y')\, s^{M-1} \int\limits_{\zeta = 0}^{s}
(1-2 \Theta \zeta + \zeta^{2})^{\lambda - \frac{1}{2}}\, d
\zeta . \label{4.2.14}
\end{equation}

For $M \geq 1$ and $\lambda \geq \frac{1}{2}$ the
integrand in \eqref{4.2.14} is continuous and $| \Theta
| \leq 1$ so $(1-2 \Theta \zeta + \zeta^{2}) \leq
(1+s)^{2}$. Therefore,
\begin{equation}
\left| K_{M} (\lambda , x, y') \right| \leq d_{1}\,K(\lam,x,y')
\,s^M\,(1+s)^{2\lam -1}.\label{4.2.14.5}
\end{equation}
The estimate
\begin{equation}
\left| K (\lambda , x, y')\right| \leq 2^{2\lam}
\sec^{2\lambda} \theta\, (| x | + | y' |
)^{-2\lambda} \label{4.2.15}
\end{equation}
is in \cite{SiegelTalvila} (Corollary 2.1).  Hence,
\begin{align}
\left| K_{M} (\lambda , x, y') \right| &\leq d_{2}\,s^M \sec^{2\lam}\theta
\,|y'|^{-2\lam}(1+s)^{-1}
\label{4.2.16}\\
&\leq d_{2}\, s^{M} \sec^{2\lambda} \theta\, | y'
|^{-2\lambda} . \label{4.2.17}
\end{align}
Multiply \eqref{4.2.16} by $| f (y')|$ and integrate
$y' \in \mathbb{R}^{n-1}$, $| y' | > 1$.  Letting $|
x | \raro \infty$, the Dominated Convergence Theorem
gives \eqref{1.10}.

When $0 < \lambda < \frac{1}{2}$ the integrand in
\eqref{4.2.14} can be singular. In this case
\begin{align}
\int\limits_{\zeta = 0}^{s} (1-2 \Theta \zeta +
\zeta^{2})^{\lambda - \frac{1}{2}} d \zeta &\leq
\int\limits_{\zeta = 0}^{s} | 1 - \zeta
|^{2\lambda -1}\, d \zeta \notag \\
&= \frac{1}{2\lambda} \begin{cases}
1-(1-s)^{2\lambda} , & 0 \leq s \leq 1 \\
1+(s-1)^{2\lambda}, &s \geq 1 \end{cases} \notag \\
&\leq \left( \frac{1}{\lambda} \right) \min (s,
s^{2\lambda}). \label{4.2.18}
\end{align}
And,
\begin{equation}
\left| K_{M} (\lambda , x, y') \right| \leq d_{3}\, s^{M}
\sec^{2\lambda} \theta\, (| x | + | y' |
)^{-2\lambda} \label{4.2.19}
\end{equation}
so \eqref{1.10} holds for $0 < \lambda < \frac{1}{2}$ as
well.

Finally, integrating \eqref{4.2.15} shows \eqref{1.10}
also holds when $M=0$.

\smallskip
\noindent
{\bf Step II}\,\,  We now prove sharpness.  Given any sequence $\{
x^{(i)}\}$ in $\Pi_{+}$ with $| x^{(i)} | \raro
\infty$ and any function $\psi (x) = o \left( | x
|^{M} \sec^{2\lambda} \theta \right)$ we find a
continuous function $f$ satisfying \eqref{1.5} for which
$F_{\lambda , M}[f] (x^{(i)})/\psi (x^{(i)}) \not \to 0$ as $i \raro
\infty$.

Note that \eqref{2.11} may be written
\begin{equation}
K_{M} (\lambda , x, y') = K(\lambda , x, y')\, s^{M}
\int\limits_{\zeta =0}^{1} (1-2 \Theta s \zeta + s^{2}
\zeta^{2})^{\lambda -1}\, \Phi_{-}(\Theta , s \zeta )\,
\zeta^{M-1}\, d \zeta . \label{4.2.20}
\end{equation}
Suppose first that $\{ x^{(i)}\}$ has a subsequence
$\tilde{x}^{(i)} = a_{i} \hat{e}_{1} + b_{i}
\hat{e}_{n}$, $i \geq 1$, where $b_{i} >0$ and $0 \leq
a_{i} \leq b_{i} \tan \theta_{0}$ for some $0 \leq
\theta_{0} < \pi /2$.  Then $0 \leq \sin \theta =
a_{i}/\sqrt{a_{i}^{2} + b_{i}^{2}} \leq \sin \theta_{0}
< 1$.  Since $\psi (x) = o(| x |^{M}
\sec^{2\lambda} \theta )$ and $1 \leq \sec \theta =
\sqrt{a_{i}^{2} + b_{i}^{2}}/b_{i} \leq \sec \theta_{0}
< \infty$ we also have $\psi (x) = o(| x |^{M})$.
We may assume that $\tilde{x}^{(i)}$ have been chosen so
that $\psi (\tilde{x}^{(i)}) \leq |
\tilde{x}^{(i)}|^{M} /i^{2}$, $i \geq 1$.

Now find a region $\Omega_{1} \subset \mathbb{R}^{n-1}$ in
which $\Phi_{-}(\Theta , s \zeta )$ is of one sign.
Consider $n \geq 3$ and $M \geq 1$.  Let $\beta_{1}$ be
the smallest positive root of $\left\{ C_{M}^{\lambda},
C_{M-1}^{\lambda} \right\}$.  And,
$C_{m}^{\lambda} (\Theta )$ is a polynomial in $\Theta$
of degree $m$ with $m$ simple zeroes in $(-1,1)$.
If $M=1$, take $\beta_1 =1$.
So 
$0<\beta_1 \leq 1$.
Now, $C_{m}^{\lambda}$ is
even or odd about the origin according as $m$ is even or
odd (\cite{Szego}, 4.7.4) and $(-1)^{m}
C_{2m}^{\lambda}(0) > 0$ (\cite{Erdelyi}, 10.9.19).
Hence, one of $\left\{ C_{M}^{\lambda},
C_{M-1}^{\lambda} \right\}$ changes sign at the origin and 
the other is of one sign on $(-\beta_1, \beta_1)$.
Therefore, for any $0 \leq \theta \leq \pi /2$,
$C_{M}^{\lambda} (\sin \theta \cos \theta' )$ and
$C_{M-1}^{\lambda} (\sin \theta\cos \theta' )$ are each
of one sign for $\arccos (\beta_{1}) \leq \theta' \leq
\pi /2$.  The same can be said when $\pi /2 \leq \theta' \leq \pi - \arccos
(\beta_{1})$.  Write $M=2 \mu + \veps_{0}$ where
$\veps_{0}$ is $0$ or $1$.  From \eqref{2.3} we see that if
$0 < t < \beta_{1}$ then ${\rm sgn} (C_{2\mu + 1}^{\lambda}
(t)) = {\rm sgn} (C_{2\mu}^{\lambda +1} (t)) = (-1)^{\mu}$
and if $- \beta_{1} < t < 0$ then ${\rm sgn} (C_{2\mu
+1}^{\lambda} (t)) = - {\rm sgn} (C_{2\mu}^{\lambda +1} (t)) =
(-1)^{\mu +1}$.  Let
\begin{eqnarray}
\Omega_{1} (\hat{y}) & = & \Bigl\{ y' \in \mathbb{R}^{n-1}
\Bigr| \arccos (\beta_{1}/2) \leq \theta' \leq \arccos (\beta_{1}/3)
\quad\text{if }  M  \text{ is even} \nonumber \\
& & \mbox{and} \quad \arccos (\beta_{1}/3) \leq \theta' \leq \pi - \arccos
(\beta_{1}/2)
\quad\text{if }  M  \text{ is odd} \Bigr\}.
\label{4.2.21}
\end{eqnarray}
Then, since $C_{M}^{\lambda}$ and $C_{M-1}^{\lambda}$
have no common roots, there exists a positive constant
$d_{4}$ such that
\begin{equation}
(-1)^{\mu + \veps_{0}} \Phi_{-}(\Theta , s \zeta ) \geq
d_{4} , \label{4.2.22}
\end{equation}
whenever $0 \leq \theta \leq \pi /2$, $y' \in \Omega_{1}
(\hat{y})$, $0 \leq \zeta \leq 1$, $s \geq 0$.  In \eqref{4.2.21}, 
$\theta'$ is restricted to lie in a smaller region than $\arccos \beta_1
\leq \theta' \leq \pi/2$ so that $(-1)^{\mu + \veps_{0}} \Phi_{-}$
will be strictly positive for $y'\in \Omega_{1}$.

From \eqref{4.2.20} we will need the estimate,
\begin{eqnarray}
(1-2 \Theta s \zeta + s^{2} \zeta^{2})^{\lambda -1} & \geq &
\left\{
\begin{array}{ll}
(1+s)^{2(\lambda -1)}, & 0 \leq \lambda \leq 1 \\
((s\zeta -\sin\theta_{0})^{2} + \cos^{2}\theta_{0})^{\lambda -1}, & 
\lambda \geq 1
\end{array}\right.
\notag \\
& \geq & (1+s)^{-2}\, \cos^{2 | \lambda -1 |} \theta_{0}.
\label{4.2.23}
\end{eqnarray}

These give
\begin{equation}
(-1)^{\mu + \veps_{0}} K_{M} (\lambda , x, y') \geq
\frac{d_{5}\, K(\lambda , x, y')\,
s^{M}}{(1+s)^{2}}, \label{4.2.24}
\end{equation}
whenever $y'\in\Omega_{1}$.

If $M=0$ then \eqref{4.2.23} and \eqref{4.2.24} hold and we
can take $\Omega_{1} = \mathbb{R}^{n-1}$.

Let
\begin{equation}
f(y') = \begin{cases}
(-1)^{\mu +\veps_{0}} f_{i} \Bigl[ 1- | y'-
 c_{i} \hat{e}_{2} | \Bigr] | y_{1}' | ; &
y' \in B_{1} (c_{i} \hat{e}_{2}), (-1)^{M} y_{1}' \geq 0
\\
0, & \text{otherwise}, \end{cases} \label{4.2.25}
\end{equation}
where $c_{i} := | \tilde{x}^{(i)} | =
\sqrt{a_{i}^{2} + b_{i}^{2}}$ and the constants $f_{i}$ are defined in
\eqref{4.2.27} below.  Then $f \!: \mathbb{R}^{n-1}  \raro
\mathbb{R}$, has support in a sequence of half balls along
the $\hat{e}_{2}$ axis and is continuous.  The factor
$[1- | y' - c_{i} \hat{e}_{2} | ] | y_{1}'
|$ makes $f$ vanish on the perimeter of the $i^{\rm
th}$ half ball.  Without loss of generality we may
assume $c_{i} \raro \infty$ monotonically so that the
$B_{1} (c_{i} \hat{e}_{2})$ are disjoint, ${\rm  supp}(f)
\subset \Omega_{1}$ and $c_{i} \geq 2$ (otherwise, take
an appropriate subsequence of $\{ \tilde{x}^{(i)} \}$).

Now, for any $j \geq 1$,
$$
F_{\lambda , M} [f] (\tilde{x}^{(j)}) \geq
\int\limits_{B_{1} (c_{j} \hat{e}_{2})} f(y') K_{M}
(\lambda , \tilde{x}^{(j)}, y')\, dy'.
$$
When $y' \in B_{1} (c_{j} \hat{e}_{2})$ we have $s =
| \tilde{x}^{(j)} | \, / | y' | \leq c_{j} /
(c_{j} -1) \leq 2$ and \newline $s \geq c_{j} /(c_{j}+1) \geq
2/3$.  And,
\begin{equation*}
\begin{split}
K(\lambda , \tilde{x}^{(j)}, y') &\geq \Bigl[ ( | y'
| + a_{j})^{2} + b_{j}^{2} \Bigr]^{-\lambda} \\
& \geq \Bigl[ (c_{j} + 1 + a_{j})^{2} + b_{j}^{2}
\Bigr]^{-\lambda} \\
& \geq (7 c_{j}^{2})^{-\lambda}.
\end{split}
\end{equation*}
Thus, using \eqref{4.2.24},
\begin{align}
F_{\lambda, M}[f] (\tilde{x}^{(j)}) & \geq \frac{d_{6}\,
f_{j}}{| \tilde{x}^{(j)}|^{2\lambda}}
\int\limits_{B_{1}} (1-| y' | )\, | y'_{1} |\, d y'
\notag \\
&= d_{7}^{-1}\, f_{j}\, | \tilde{x}^{(j)}|^{-2 \lambda} .
\label{4.2.26}
\end{align}
Let
\begin{equation}
f_{i} = d_{7}\, \psi (\tilde{x}^{(i)})\, |
\tilde{x}^{(i)} |^{2 \lambda} . \label{4.2.27}
\end{equation}
Then $f_{i}/c_{i}^{M+2\lambda} \leq (d_{7} i^{2})^{-1}$
and ${\dss\sum_{i=1}^{\infty}} f_{i} c_{i}^{-(M+2\lambda
)} < \infty$ so \eqref{1.5} holds.  And, on $\{ \tilde{x}^{(i)} \}$ we 
have $F_{\lambda , M}[f](\tilde{x}^{(j)}) \geq \psi (\tilde{x}^{(j)})$ for 
each $j \geq 1$ so ${\dss\operatornamewithlimits{lim \, sup}_{i \raro \infty}}\,
F_{\lambda ,M}[f](x^{(i)})/\psi
(x^{(i)}) \geq 1$ and 
$F_{\lambda , M} [f] (x^{(i)})/ \psi (x^{(i)}) \not\to 0$ as $i \raro \infty$.

When $n =2$, write $x_{1} = r \cos \phi$, $x_{2} = r \sin
\phi$.  Then in place of \eqref{2.11}  we have\newline 
$F_{\lambda , M}[f] (x) =
{\int_{-\infty}^{\infty}} f(\xi )$ $K_{M} (\lambda
, x, \xi )\,d \xi$ where
$$
K_{M} (\lambda , x, \xi ) = K (\lambda , x, \xi )\left(
\frac{r}{\xi} \right)^{M} \!\int\limits_{\zeta = 0}^{1}\!
\left( 1-2 \frac{r\zeta}{\xi} \cos \phi +
\frac{r^{2}\zeta^{2}}{\xi^{2}} \right)^{\lambda -1} \!\Phi_{-}\!
\left(\cos \phi , \frac{r\zeta}{\xi} \right)
\zeta^{M-1}\, d \zeta .
$$
If $0 \leq \theta \leq \theta_{0} < \pi /2$ then $0 <
\phi_{0} \leq \phi \leq \pi - \phi_{0} < \pi$ where
$\phi_{0} = \pi /2 - \theta_{0}$.

Let $t_{i}$, $1 \leq i \leq q$, be the roots of
$C_{M}^{\lambda} \circ \cos$ and $C_{M-1}^{\lambda}
\circ \cos$ in $[\phi_{0},\, \pi - \phi_{0}]$, ordered by
size.  We then have the partition $\phi_{0} = t_{0} \leq
t_{1} < t_{2} < \!\cdots\!$ $< t_{q-1} < t_{q} \leq t_{q+1}
= \pi - \phi_{0}$.  In each interval $[t_{i}, t_{i+1}]$,
$0 \leq i \leq q$, $C_{M}^{\lambda} \circ \cos$ and
$C_{M-1}^{\lambda} \circ \cos$ are each of one sign. If
$\phi_{0}$ is a root, we omit the singleton $\{ t_{1}
\}$, similarly with $\pi - \phi_{0}$.

For any sequence $\phi_{i} \in [ \phi_{0},\, \pi - \phi_{0}
]$, $i \geq 1$, there is a subsequence $\{
\tilde{\phi}_{i} \}$ in one of the above intervals $[
t_{j}, t_{j+1}]$.  If $C_{M}^{\lambda} (\cos
\tilde{\phi}_{i})$ and $C_{M-1}^{\lambda} (\cos
\tilde{\phi}_{i})$ are of the same sign, take $\Omega_{1} =
\{ \xi \in \mathbb{R} | \xi < 0 \}$ and $\Omega_{1} =
\{ \xi \in \mathbb{R} | \xi > 0 \}$ if they are of
opposite sign.  Then $(-1)^{\mu_{0}} \Phi_{-}(\cos
\tilde{\phi}_{i} , r \zeta / \xi ) \geq 0$ for $i \geq
1$, $\xi \in \Omega_{1}$, where $(-1)^{\mu_{0}} = {\rm sgn}
(C_{M}^{\lambda} (\cos \tilde{\phi}_{i}))$ $(\mu_{0} = 0$
or $1$).  Since $C_{M}^{\lambda}$ and
$C_{M-1}^{\lambda}$ have no common zeroes there is a
subsequence $\{ \check{\phi}_{i} \}$ of $\{
\tilde{\phi}_{i} \}$ such that either $C_{M}^{\lambda}
(\cos \check{\phi}_{i})$ or $C_{M-1}^{\lambda}
(\cos \check{\phi}_{i})$ is bounded away from
zero for all $i \geq 1$.  Hence, there is a positive
constant $\check{d}_{5}$ such that $(-1)^{\mu_{0}} \Phi_{-}
(\cos \check{\phi}_{i}, \check{r}_{i}
\zeta / \xi ) \geq \check{d}_{5} (\check{r}_{i}
\zeta / | \xi | )^{\mu_{1}}$ for $i \geq 1$.  Here
$\check{x}^{(i)} = \check{r}_{i} \cos
\check{\phi_{i}}\, \hat{e}_{1} + \check{r}_{i}
\sin \check{\phi}_{i}\, \hat{e}_{2}$ is a
sub-subsequence of the given sequence $\{ x^{(i)} \}$ and
$\mu_{1}$ is $0$ or $1$.
We now proceed in a similar manner to the case $n \geq
3$ given above.

\smallskip
\noindent
{\bf Step III}\,\,  In the previous argument $0 \leq \theta_{0} < \pi /2$
was arbitrary so now suppose that given the sequence $\{
x^{(i)} \}$ there is a subsequence $\tilde{x}^{(i)} =
a_{i} \hat{e}_{1} + b_{i} \hat{e}_{n}$ such that $\sin
\theta_{0} \leq \sin \theta = a_{i} /\sqrt{a_{i}^{2} +
b_{i}^{2}} <1$.  Since $0 < b_{i} \leq a_{i} \cot
\theta_{0}$ we may assume $0 < b_{i} \leq a_{i}/2$ and
that $a_{i} \raro \infty$ monotonically.

Find a region $\Omega_{2} \subset \mathbb{R}^{n-1}$ on which
$K_{M}$ is of one sign.  Let $M \geq 1$ and let $x =
x_{1} \hat{e}_{1} + x_{n} \hat{e}_{n} \in \{
\tilde{x}^{(i)} \}$.  We will take
$1 < A < 2$ close enough to $1$ so that if
\begin{equation}
\frac{1}{A} \leq s \leq A, \quad \frac{1}{A} \leq \Theta
\leq 1 \label{4.2.28}
\end{equation}
then $K_{M}
(\lambda , x, y')$, $C_{M}^{\lambda} (\Theta )$ and
$C_{M-1}^{\lambda} (\Theta )$ are positive.  From
\eqref{1.1} and \eqref{1.6}, we have $K_{M} (\lambda , x, y')
\geq | y' |^{-2\lambda} \Bigl[ (1-2 A^{-2} +
A^{2})^{-\lambda} - \gamma_{\lambda , M}^{-\lambda}
\Bigr]$, where $\gamma_{\lambda , M} = \left(
{\dss\sum_{m=0}^{M-1}} 2^{m} C_{m}^{\lambda}(1)
\right)^{-1/\lambda}$ and $0 < \gamma_{\lambda , M} <
\infty $.  Note that $A >1$ implies $1-2 A^{-2} +
A^{2} > 0$.  Now, $K_{M} > 0$ if $A^{4} +
(1-\gamma_{\lambda , M}) A^{2} - 2 < 0$.  Let $r_{0} > 1$
be the largest root of this quartic.
Let
$\beta_{2}$ be the largest zero of $C_{M}^{\min(1,\lambda)}$.  
Then $\cos(\pi/(M+1)) \leq \beta_{2}
\leq \cos (\pi /(2M))$ (\cite{Szego}, 6.21.7).  Hence,
if $1 < A < \min (2, r_{0}, \sec (\pi / (2M))$ and $s$
and $\Theta$ are as in \eqref{4.2.28} then $K_{M} (\lambda
, x, y') > 0$, $C_{M}^{\lambda} (\Theta ) > 0$ and
$C_{M-1}^{\lambda} (\Theta ) > 0$ (\cite{Szego}, 6.21.3).

To satisfy $\Theta \geq 1/A$ in \eqref{4.2.28}, we will
restrict $x$ and $y'$ so that $\sin \theta \geq
1/\sqrt{A}$ and $\cos \theta' \geq 1/\sqrt{A}$.  First,
take $\theta_{0} = \arcsin \left( \sqrt{A/(2A-1)}
\right)$ then $\sin \theta \geq \sin \theta_{0} =
\sqrt{A/(2A-1)} \geq 1/\sqrt{A}$.  And, since $y = x_{1}
\hat{e}_{1}$, we have $\cos \theta' = (y\!\cdot\! y')/(| y
| \, | y' | ) = {\hat e}_{1}\!\cdot\!y'\, | y' |^{-1}$ for $y'
\neq 0$.  Let
\begin{equation}
\Omega_{2} (\hat{y}) = \left\{ y' \in \mathbb{R}^{n-1} \bigl|
| y' | > 1, \,\, 1/\sqrt{A} < \cos \theta' \leq 1
\right\} , \label{4.2.29}
\end{equation}
a portion of a cone with axis along $\hat{e}_{1}$.
If $y' \in \Omega_{2}$ then $\cos \theta' \geq
1/\sqrt{A}$.  If $n=2$, take $\Omega_{2} = \{ \xi \in
\mathbb{R} | \xi > 1 \}$. 

Define $f\!: \mathbb{R}^{n-1} \raro \mathbb{R}$ by
\begin{equation}
f(y') = \begin{cases}
&f_{i} \left
( 1 - \frac{1}{b_{i}} | y' - a_{i} \hat{e}_{1} |
\right), \, y' \in B_{b_{i}} (a_{i} \hat{e}_{1}) \quad
\text{for some} \quad i \geq 1 \\
&(-1)^{M} A_{\lambda} f_{i} \left( 1 - \frac{1}{b_{i}}
| y' +a_{i} \hat{e}_{1} | \right), \, y' \in
B_{b_{i}} (-a_{i} \hat{e}_{1}) \quad \text{for some}
\quad i \geq 1 \\
& 0, \quad \text{otherwise}, \end{cases} \label{4.2.30}
\end{equation}
where $A_{\lambda} \geq 1$ is given in \eqref{4.2.32} and
$f_{i}$ in \eqref{4.2.38}.  By taking an appropriate
subsequence of $\{ \tilde{x}^{(i)} \}$ we may
assume the balls $B_{b_{i}} (a_{i} \hat{e}_{1})$ are
disjoint $(a_{i+1} \geq 3 a_{i}$ suffices).  The
condition $\sin \theta_{0} = \sqrt{A/(2A-1)}$ ensures
that each $B_{b_{i}} (a_{i} \hat{e}_{1}) \subset
\Omega_{2}$.  Then $f$ is continuous, has support on a
sequence of balls along the $\hat{e}_{1}$ axis and is
non-negative for $y_{1}' \geq 0$.

With $y' \in \Omega_{2}$ such that $y_{1}' >0$ and $x$
as above (preceding \eqref{4.2.28}), $A^{-1} \leq \Theta =
\sin \theta \cos \theta' \leq 1$.  So $C_{M}^{\lambda}
(\Theta )$, $C_{M-1}^{\lambda} (\Theta ) > 0$.  As a
function of $s$, with fixed $\Theta$ as in \eqref{4.2.28},
the integral in \eqref{2.11} is zero when $s=0$,
is an increasing function of $s$ for
 $0 < s < M C_{M}^{\lambda} (\Theta ) / [(2
\lambda + M-1) C_{M-1}^{\lambda} (\Theta )]$ (where it
has a maximum) and decreases for larger values of $s$.
And, we know from the analysis following \eqref{4.2.28}
that this integral is positive at $s=A$. Hence, $K_{M}
(\lambda , x, y') >0$ for $0 < s \leq A$ (with $y' \in
\Omega_{2}$,\,\,$y_{1}' > 0$).

If $y' \in \Omega_{2}$ and $y_{1}' < 0$ then $\Theta <
0$.  Since $C_{m}^{\lambda} (-t) = (-1)^{m}
C_{m}^{\lambda} (t)$ we have $\Phi_{-} (\Theta , \zeta ) =
(-1)^{M} \Phi_{+} (| \Theta | , \zeta )$ and ${\rm sgn}
(K_{M} (\lambda , x, y')) = (-1)^{M}$.  From
\eqref{4.2.30},\newline $f(y') K_{M}$$(\lambda , x, y') \geq 0$.

Define
\begin{equation}
\Omega_{\gtrless} = \{ y' \in \Omega_{2} | y_{1}'
\gtrless 0, \quad s > A \} . \label{4.2.31}
\end{equation}
If $x \in \{ \tilde{x}^{(i)} \}$ then $f(y') K_{M}
(\lambda , x, y') \geq 0$ for $y' \in \Omega_{2}$ except
possibly for $y' \in \Omega_{>}$.  By taking
$A_{\lambda}$ large enough we can ensure
${\dss\int\limits_{\Omega_{>} \cup\, \Omega_{<}}} f(y')
K_{M} (\lambda , x, y')\, dy' \geq 0$.  Indeed, let
$y^{*}$ be the reflection of $y'$ in the hyperplane
$y_{1}' =0$ and $\theta^{*}$ the angle between $y^{*}$
and $y$. Then $y^{*} \in \Omega_{<}$ if and only if $y'
\in \Omega_{>}$.

If $\lambda \geq 1$ and $y^{*} \in
\Omega_{<}$ then, as in \eqref{2.10}, $\Theta^{*} := \sin
\theta \cos \theta^{*} = - \Theta$.  Then, using
\eqref{2.11} and \eqref{4.2.30},
\begin{align*}
f(y^{*}) K_{M} (\lambda , x, y^{*}) &= A_{\lambda} f(y')
\Bigl[ | y' |^{2} + 2 \Theta | y' | \, |
x | +| x |^{2} \Bigr]^{-\lambda}
\int\limits_{\zeta = 0}^{s} \frac{\Phi_{+} (\Theta ,
\zeta )\, \zeta^{M-1}\, d \zeta}{(1+2 \Theta \zeta +
\zeta^{2})^{1-\lambda}} \\
&\geq A_{\lambda} f(y') | x |^{-2\lambda}
(1+A^{-1})^{-2\lambda} \int\limits_{\zeta =0}^{s}
(1+\zeta^{2})^{\lambda - 1} \Phi_{+} (\Theta , \zeta )\,
\zeta^{M-1}\, d \zeta .
\end{align*}
And,
\begin{align*}
f(y') K_{M} (\lambda, x, y') &= f(y') \Bigl[ | y'
|^{2} -2 \Theta | y' | \, | x | + | x
|^{2} \Bigl]^{-\lambda} \int\limits_{\zeta =0}^{s}
\frac{\Phi_{-} (\Theta , \zeta )\, \zeta^{M-1}\, d
\zeta}{(1-2 \Theta \zeta + \zeta^{2})^{1-\lambda}} \\
&\geq -f (y') | x |^{-2\lambda}
(1-A^{-1})^{-2\lambda} \int\limits_{\zeta =0}^{s}
(1+\zeta^{2})^{\lambda -1}\, \Phi_{+} (\Theta , \zeta )\,
\zeta^{M-1}\, d \zeta .
\end{align*}
Therefore,
$$
\int\limits_{\Omega_{<} \cup\, \Omega_{>}} f(y') K_{M}
(\lambda , x, y')\, dy' \geq 0 \quad \text{if} \quad
A_{\lambda} \geq (A+1)^{2\lambda} (A-1)^{-2\lambda}.
$$
If $0 < \lambda < 1$ and $y^{*} \in \Omega_{<}$ then
\begin{align*}
f(y^{*}) K_{M} (\lambda , x, y^{*}) &\geq A_{\lambda}\,
f(y')\, (| x | + | y' | )^{-2 \lambda}\,
(1+s)^{2\lambda -2} \int\limits_{\zeta = 0}^{s} \Phi_{+}
(\Theta , \zeta )\, \zeta^{M-1}\, d \zeta \\
&\geq A_{\lambda}\, f(y')\, | x |^{-2\lambda} \left( 1
+ A^{-1} \right)^{-2 \lambda} s^{M+2 \lambda - 1}\,
\frac{(2 \lambda + M-1)\, C_{M-1}^{\lambda}(\Theta )}{4(M+1)}.
\end{align*}
If $0< \lambda <1/2$ then, using \eqref{4.2.13} and \eqref{4.2.18},
$$
f(y') K_{M} (\lambda , x, y') \geq -f(y') (| x |
- | y' | )^{-2 \lambda}\, 2 \lambda \binom{
2 \lambda + M}{M-1}\, \frac{s^{M + 2\lambda-1}}{2\lambda}.
$$
And, if $1/2 \leq \lambda <1$, 
$$ 
f(y') K_{M} (\lambda , x, y') \geq -f(y') (| x |
- | y' | )^{-2 \lambda}\, 2 \lambda \binom{
2 \lambda + M}{M-1}\,
\frac{s^{M}}{M}\left(1+s^{2}\right)^{\lambda - \frac{1}{2}}.
$$
Hence, for $0 < \lambda < 1$,
$$
f(y') K_{M} (\lambda , x, y')
\geq -2{\sqrt 2}\,f (y') | x |^{-2\lambda} \left( 1 -
A^{-1} \right)^{-2\lambda} s^{M+2 \lambda -1}
\binom{2 \lambda +M}{M-1}.
$$
And, ${\dss\int\limits_{\Omega_{<} \cup\, \Omega_{>}}}
f(y') K_{M} (\lambda , x, y')\, d y' \geq 0$ if
$$
A_{\lambda} \geq \frac{8{\sqrt 2}\,(M+1)}{2 \lambda +M-1}
\binom{2 \lambda +M }{M-1} \left(
\frac{A+1}{A-1} \right)^{2 \lambda} \left[
\min_{A^{-1} \leq t \leq 1} C_{M-1}^{\lambda} (t)
\right]^{-1} .
$$

Hence, for $\lambda > 0$, $x \in \{ \tilde{x}^{(i)} \}$, if we take
\begin{equation}
A_{\lambda} \geq \left( \frac{A+1}{A-1} \right)^{2
\lambda} \max \left( 1, \frac{8{\sqrt 2}\,(M+1)}{2 \lambda +M-1}
\binom{2 \lambda +M }{M-1} \left[
\min_{A^{-1} \leq t \leq 1} C_{M-1}^{\lambda} (t)
\right]^{-1} \right) \label{4.2.32}
\end{equation}
then $F_{\lambda ,M} [f] (x) \geq
{\dss\int\limits_{\Omega_{3}}} f(y') K_{M} (\lambda , x,
y')\, dy'$, where
\begin{equation}
\Omega_{3} = \bigl\{ y' \in \Omega_{2} \bigl| y_{1}' >
0, \, A^{-1} < s < A \bigr\} . \label{4.2.33}
\end{equation}
Note that if $x = a_{i} \hat{e}_{1} + b_{i} \hat{e}_{n}$
then $B_{b_{i}} (a_{i} \hat{e}_{1}) \subset \Omega_{3}$
if $a_{i} - | x| /A \geq b_{i}$ and $A| x |
- a_{i} \geq b_{i}$. Since $a_{i} = | x | \sin
\theta$, $b_{i} = | x | \cos \theta$ and
$\theta_{0} \leq \theta < \pi /2$, these conditions are
satisfied if $\frac{1}{2} [ \pi - \arcsin (1-A^{-2})] \leq
\theta_{0} < \pi /2$, i.e., by taking $\theta_{0}$ close
enough to $\pi /2$.

From \eqref{2.11},
\begin{equation}
K_{M} (\lambda , x, y') = K(\lambda , x, y')\, s^{M}
\int\limits_{\zeta =0}^{1} (1-2 \Theta s \zeta + s^{2}
\zeta^{2})^{\lambda -1}\, \Phi_{-} (\Theta , s \zeta )\,
\zeta^{M-1}\, d \zeta, \label{4.2.34}
\end{equation}
which is strictly positive on $\ol{\Omega}_{3}$
(\eqref{4.2.28} and following).  And, $K(\lambda , x, y')$
is positive but singular at $s = \Theta = 1$.  Using
\eqref{2.2}, the integral \eqref{4.2.34} above reduces to
$$
\frac{\Gamma(2 \lambda +M)}{\Gamma(2\lambda)\, \Gamma(M)}
\int\limits_{\zeta = 0}^{1} (1-\zeta)^{2\lambda -1}
\zeta^{M-1}\, d \zeta > 0
$$
at $s = \Theta = 1$.  The integral in \eqref{4.2.34} is a
strictly positive continuous function of $s$ and
$\Theta$ when the conditions in \eqref{4.2.28} are
satisfied.  Hence, it must be bounded below by a
positive constant, say $d_{8}$, i.e.,
\begin{equation}
K_{M} (\lambda , x, y') \geq d_{8}\, K(\lambda , x, y')\,
s^{M} \quad \text{for} \quad y' \in \Omega_{3}, \, \sin
\theta \geq \sin \theta_{0} . \label{4.2.35}
\end{equation}

If $M =0$ we can dispense with the sets $\Omega_{2},
\Omega_{3}, \Omega_{<}$ and $\Omega_{>}$.  In
\eqref{4.2.30}, $A_{\lambda} = 1$ and $f$ is extended as
an even function.  Then \eqref{4.2.35} holds for $x \in
\Pi_{+}$ with $d_{8} =1$.

For $M \geq 0$ each element of the sequence
$\tilde{x}^{(i)} = a_{i} \hat{e}_{1} + b_{i}
\hat{e}_{n}$ satisfies $\sin \theta \geq \sin
\theta_{0}$ so, using \eqref{4.2.30} and \eqref{4.2.35}
\begin{align}
F_{\lambda , M} [f] (\tilde{x}^{(j)}) & \geq d_{8} |
\tilde{x}^{(j)}|^{M} \int\limits_{\Omega_{3}} f(y') K
(\lambda , \tilde{x}^{(j)}, y')\, | y' |^{-M}\,
d y' \notag \\
&\geq d_{8} (a_{j}^{2} + b_{j}^{2} )^{M/2} f_{j}
\int\limits_{B_{b_{j}}(a_{j} \hat{e}_{1})} \frac{( 1 -
| y' - a_{j} \hat{e}_{1} |\,b_{j}^{-1})}{
(| y' - a_{j} \hat{e}_{1} |^{2} +
b_{j}^{2})^{\lambda}}\,
| y' |^{-M}\, dy' \notag\\
&\geq \frac{d_{8} (a_{j}^{2} + b_{j}^{2})^{M/2}\,
f_{j}\,b_{j}^{n-1}}{(a_{j} + b_{j})^{M}\, b_{j}^{2\lambda}}
\int\limits_{B_{1}} (1- | y' | )\, (| y' |^{2}
+1)^{-\lambda}\, d y' \notag \\
& \geq d_{9}\, f_{j}\, b_{j}^{n-1-2\lambda} \label{4.2.36}
\end{align}
where $d_{9} = d_{8}\, 2^{-M/2} (n-1) \om_{n-1}
{\dss\int\limits_{\rho=0}^{1}} (1-\rho) (\rho^{2}+1)^{-\lambda}\,
\rho^{n-2}\, d\rho$.

Note that \eqref{1.5} holds if and only if
\begin{equation}
\sum_{i=1}^{\infty} \frac{f_{i} b_{i}^{n-1}}{a_{i}^{M+2
\lambda}} < \infty . \label{4.2.37}
\end{equation}

Now suppose $\psi \!: \mathbb{R}^{n-1} \raro (0,\infty)$ such that
$\psi (x) = o (| x |^{M} \sec^{2 \lambda} \theta
)$.  On the sequence $\tilde{x}^{(j)} = a_{j}
\hat{e}_{1} + b_{j} \hat{e}_{n}$, $| x |^{M}
\sec^{2 \lambda} \theta = | x |^{M+2\lambda}
x_{n}^{-2\lambda} = (a_{j}^{2}+ b_{j}^{2})^{M/2+\lambda}\,
b_{j}^{-2\lambda}$ and $\psi (\tilde{x}^{(j)}) = o
(a_{j}^{M + 2\lambda}\,b_{j}^{-2\lambda})
$\, (since $0 < b_{i} \leq a_{i}/2$).
We may assume that $\{ \tilde{x}^{(i)} \}$ has been
chosen so that $\psi (\tilde{x}^{(i)}) \leq
a_{i}^{M+2\lambda}\, b_{i}^{-2\lambda}\, i^{-2}$ for $i \geq
1$.   Let
\begin{equation}
f_{i} = d_{9}^{-1}\, \psi (a_{i} \hat{e}_{1} + b_{i}
\hat{e}_{n})\, b_{i}^{2 \lambda - n+1}, \quad i \geq 1 .
\label{4.2.38}
\end{equation}
Then \eqref{4.2.37} is satisfied and $F_{\lambda ,M} [f]
(\tilde{x}^{(j)}) \geq \psi (\tilde{x}^{(j)})$ so
$F_{\lambda , M}[f] (x^{(i)}) / \psi (x^{(i)}) \not \to 0$ as 
$i \raro\IN$ and 
$F_{\lambda , M}[f] (x) \neq o(|x|^{M}\sec^{2\lambda}\theta)$.
Hence, the order relation $F_{\lambda , M}[f] (x) = o
(| x |^{M} \sec^{2\lambda} \theta )$ is sharp for
$x \in \Pi_{+}$ of form $x = x_{1} \hat{e}_{1} + x_{n}
\hat{e}_{n}$.

\smallskip
\noindent
{\bf Step IV}\,\,
For $n =2$ this completes the proof.  For $n \geq 3$ we
now remove this restriction on $x$.  For any sequence
$\{ x^{(i)} \}$ in $\Pi_{+}$, we can write $x^{(i)} = |
x^{(i)} |\, \hat{x}^{(i)}$ where $\hat{x}^{(i)} \in
\partial B_{1}^{+} ={\{ x \in \mathbb{R}^{n} \bigl| | x
| = 1, \, x_{n} > 0 \}}$. Then $\{ \hat{x}^{(i)} \}$
must have a limit point, say $\hat{s}_{0}$, in the
compact set $\ol{\partial B_{1}^{+}}$.  Let $\theta_{0}$ be
the angle between $\hat{s}_{0}$ and $\hat{e}_{n}$.

If $0 < \theta_{0} < \pi /2$ then let $\hat{s}_{1}$ be in
the direction of the projection of $\hat{s}_{0}$ onto
$\partial \Pi_{+}$ and let $\hat{s}_{2} \in \partial
\Pi_{+}$ be any unit vector orthogonal to
$\hat{s}_{1}$.  For any $\delta > 0$ there is a
subsequence $\tilde{x}^{(i)}_{\delta} =\tilde{x}^{(i)} +
\delta_{i} \hat{t}_{i}$ where $\tilde{x}^{(i)} = a_{i}
\hat{s}_{1} + b_{i} \hat{e}_{n}$, $c_{i} = |
\tilde{x}^{(i)} | \raro \infty$ monotonically, $b_{i}
> 0$, $\{ \hat{s}_{1}, \hat{e}_{n}, \hat{t}_{i} \}$ is
orthonormal, each $\hat{t}_{i} \in \partial \Pi_{+}$ and
$0 \leq \delta_{i} \leq \delta$. 
We can now try to repeat the first part
of the sharpness proof, beginning with \eqref{4.2.20}.
Then $\hat{s}_{1}$ and $\hat{s}_{2}$ play the roles
$\hat{e}_{1}$ and $\hat{e}_{2}$ did before, except that
we now have perturbations by $\delta_{i}$. 

Let $x\in\{{\tilde x}_{\delta}^{(i)}\}$.  Let $\eta_{i}$ be the
angle between $a_{i}{\hat s}_{1} + \delta_{i} {\hat t}_{i}$
and ${\hat s}_{1}$.  Without loss of generality $a_{i} \geq 1$.
We have $0 \leq \eta_{i} = \arctan{(\delta_{i}/a_{i})} \leq \delta$.
Hence, we can replace \eqref{4.2.21} with the narrower cone
\begin{eqnarray}
\Omega'_{1} (\hat{y}) & = & \Bigl\{ y' \in \mathbb{R}^{n-1}
\Bigr|\arccos (\beta_{1}/2) + \delta \leq \theta_{1}' \leq \arccos 
(\beta_{1}/3) - \delta
\quad \text{if } M  \text{ is even} \nonumber \\
& & \!\!\mbox{and} \,\,\arccos (\beta_{1}/3) + \delta \leq \theta_{1}' \leq \pi - \arccos
(\beta_{1}/2) - \delta
\,\,\,\, \text{if }  M  \text{ is odd} \Bigr\}
\end{eqnarray}
where $\theta'_{1}$ is the angle between $y'$ and $a_{i}{\hat s}_{1}$.
(Take $2\delta < \arccos (\beta_{1}/3) -\arccos{(\beta_{1}/2)}$.)  For any 
$x \in \{
{\tilde x}_{\delta}^{(i)}\}$, if $ y' \in \Omega'_{1}$ then $|\cos{\theta'}|
\leq \beta_{1}/2$ and \eqref{4.2.22} holds.

If ${\hat s}_{0} = {\hat e}_{n}$ ($\theta_{0} = 0$) then take
$a_{i}\equiv 0$.  Let ${\hat s}_{1} = {\hat e}_{1}$ and ${\hat s}_{2} = {\hat e}_{2}$.  For any $x\in\{{\tilde x}_{\delta}^{(i)}\}$ and $y' \in \Omega'_{1}$
we have $0 \leq \theta \leq \delta$ and so $0 \leq \sin{\theta} \leq 
\sin{\delta}$.  Therefore, $|\Theta| = |\sin{\theta}\cos{\theta'}| \leq
\delta \leq \beta_{1}/2$ for small enough $\delta$.  And, \eqref{4.2.22}
holds.

Now, for $0 \leq \theta_{0} < \pi/2$, replace \eqref{4.2.25} with
\begin{equation}
f(y') = \begin{cases}
(-1)^{\mu +\veps_{0}} f_{i} \Bigl[ 1- |y'-
 c_{i} \hat{e}_{\delta} | \Bigr] y_{\delta}'\, ; &
y' \in B_{1} (c_{i} \hat{e}_{\delta}),\,\, y'_{\delta} \geq 0
\\
0, & \text{otherwise}, \end{cases} 
\end{equation}
where $y'_{\delta} = y' \!\cdot\!((-1)^{M}\cos\delta\, {\hat s}_{1}
- \sin\delta\,{\hat s}_{2})$.  We align the half balls
of the support of $f$ along the unit vector $\hat{e}_{\delta}$ in the
direction $(-1)^{M}\sin{\delta}\,{\hat s}_{1} + \cos{\delta}\,{\hat s}_{2}$
so that $B_{1}(c_{i}\hat{e}_{\delta}) \subset \Omega'_{1}$.

The rest of the proof for this case follows without serious change, 
through \eqref{4.2.27}.

If $\hat{s}_{0} \in \partial \Pi_{+}\,\, (\theta_{0} = \pi/2)$ then
$\hat{s}_{1} = \hat{s}_{0}$ and a subsequence approaches the boundary.
As before, for any $\delta > 0$ there is a subsequence of form $\tilde{x}
_{\delta}^{(i)}={ \tilde{x}^{(i)} + \delta_{i}\hat{t}_{i}}$ where
\newline
$\tilde{x}^{(i)} = a_{i}\hat{s}_{1} + b_{i}\hat{e}_{n}, \,\, 0 < b_{i} \leq
a_{i}, a_{i}\geq 1,  a_{i} \raro \infty \mbox{ monotonically, } b_{i}/a_{i}
\raro 0, \{\hat{s}_{1},\hat{e}_{n}, \hat{t}_{i}\} \mbox{ is }\newline
\mbox{ orthonormal }\mbox{ and }
0 \leq \delta_{i} \leq \delta$.  Follow the second part of the sharpness
proof, from \eqref{4.2.28}.

Let
$$
B_{\delta} = \min\left(\frac{A}{(1 + \delta {\sqrt A})^{2}},\, A-\delta,\,
\frac{A}{1+\delta A}\right)
$$
then $B_{\delta} < A$.  And, $B_{\delta} >1$ if
$$\begin{array}{ccccl}
0 & < & \delta & < & \min\left(({\sqrt A} -1)/\sqrt A,\, A-1,\, (A-1)/A \right) \\
  &   &        & = & ({\sqrt A}-1)/{\sqrt A}.
\end{array}
$$
Without loss of generality, we can take $A$ satisfying
the conditions following \eqref{4.2.28} and $0 < \delta < (A -1)/A < 1/2$.

For each $j \geq 1$, let $\theta'$ be the angle between $y'$ and 
$a_{j}\hat{s}_{1}$
and $\theta'_{\delta}$ the angle between $y'$ and 
$a_{j}\hat{s}_{1} + \delta_{j}\hat{t}_{j}$.  We have 
$\theta' -\delta \leq \theta'_{\delta} \leq \theta' + \delta$
so replace \eqref{4.2.29} with
\begin{equation}
\Omega'_{2} (\hat{y}) = \left\{ y' \in \mathbb{R}^{n-1} \bigl|
 |y'| > 1, \,\, 0 \leq \theta'_\delta < 
\arccos(1/\sqrt {B_{\delta}})-\delta\right\}.
\end{equation}
For each $j \geq 1$, let $\theta$ be the angle between $\tilde{x}
^{(j)}$ and $ \hat{e}_{n}$ and $\theta_{\delta}$ the angle between
$\tilde{x}_{\delta}
^{(j)}$ and $ \hat{e}_{n}$.  Then
\[
\sin\theta_{\delta} = \frac{|a_{j}\hat{s}_{1} + \delta_{j}\hat{t}_{j}|}
{|\tilde{x}
^{(j)} + \delta_{j}\hat{t}_{j}|} \geq \frac{a_{j} - \delta}
{|\tilde{x}
^{(j)}| + \delta}.
\]
For large enough $|\tilde{x}^{(j)}|$, we have 
$\sin\theta_{\delta}  \geq  a_{j}/|\tilde{x}
^{(j)}| - \delta 
=  \sin\theta -\delta$.
It follows from the first component of the definition of $B_{\delta}$
that $\sin\theta \geq 1/\sqrt {B_{\delta}}$ implies $\sin\theta_{\delta}
\geq 1/\sqrt A$.

Write $s=|\tilde{x}
^{(j)}|/|y'|$, $s_{\delta} = |\tilde{x}
^{(j)}_{\delta}|/| y'|$.  Then for $y' \in \Omega'_{2}$, $s-\delta \leq
s_{\delta} \leq s+\delta$.  From the second and third components in the
definition of $B_{\delta}$, $1/B_{\delta} \leq s \leq B_{\delta}$ implies
$1/A \leq s_{\delta} \leq A$.  Hence, we can replace $\Omega_{2}$ with 
$\Omega'_{2}$
and carry out the sharpness proof for $a_{j}\hat{s}_{1} + b_{j}\hat{e}_{n}$
with the following changes.  In \eqref{4.2.30}, replace $a_{i}\hat{e}_{1}$
with $a_{i}\hat{s}_{1} + \delta_{i}\hat{t}_{i}$.  In \eqref{4.2.31} and
\eqref{4.2.33}, replace $y'_{1}$ with $y'\!\cdot\! \hat{s}_{1}$.  The rest of
the proof, through \eqref{4.2.38}, follows with minor changes.
$\qed$

The growth estimate on $\Ff$ gives estimates for the solutions of the half space
Dirichlet and Neumann problems.  The modified kernel introduces
a singularity at the origin of the integration space.  To avoid
integrating $f$ there, a continuous
cutoff function that vanishes in a neighbourhood of the origin is used.

\begin{cor}\label{cor2.1}
Let $w \!: \mathbb{R}^{n-1} \raro [0,1]$ be continuous such
that $w(y) \equiv 0$ when $| y | \leq 1$ and $w(y)
\equiv 1$ when $| y | \geq 2$.  Let $f$ be
continuous on $\mathbb{R}^{n-1}$ and satisfy \eqref{1.5} with
$\lambda = n/2\quad (n \geq 2)$.  The function $u(x) = D_{M}
[wf] (x) +D [(1-w)f](x)$ satisfies
\begin{align}
u& \in C^{2} (\Pi_{+}) \cap C^{0} (\ol{\Pi}_{+})
\label{2.39}\\
%\end{equation}
%\begin{equation}
&\Delta u = 0, \quad x \in \Pi_{+} \label{2.40}\\
%\end{equation}
%\begin{equation}
&u = f, \quad x \in \partial \Pi_{+} \label{2.41}\\
%\end{equation}
%\begin{equation}
u(x) = o \bigl( | x |^{M+1}&\sec^{n-1} \theta
\bigr); \quad x \in \Pi_{+}, \quad | x| \raro
\infty . \label{2.42}
\end{align}
\end{cor}

\bigskip
\noindent
{\it Proof}:  That $u$ is a classical solution,
\eqref{2.39}, \eqref{2.40}, \eqref{2.41}, is contained in
Corollary 2 of \cite{Yoshida}.  To prove \eqref{2.42},
note that the Theorem gives $D_{M} [wf] (x) = \alpha_{n} x_{n}
F_{\frac{n}{2}, M} [wf] (x) = o\bigl( | x
|^{M+1} \sec^{n-1} \theta \bigr)$. And,
\begin{align*}
\left| D[(1-w)f] (x) \right| & \leq \alpha_{n} x_{n}
\int\limits_{| y' | < 2} | f(y') |\, (| x
| - 2)^{-n}\, d y' \\
& \leq \alpha_{n}\, 2^{n} x_{n} | x |^{-n} \quad
\text{if} \quad | x | \geq 4,
\end{align*}
so \eqref{2.42} is satisfied. \qed

\begin{cor}\label{cor2.2}
Let $f$ and $w$ be as in Corollary \ref{cor2.1} such
that \eqref{1.5} holds with\newline ${\lambda = (n-2)/2}\quad(n \geq
3)$. Then $v(x) = N_{M} [wf] (x) + N [(1-w)f](x)$
satisfies \eqref{2.40} and
\begin{align}
v & \in C^{2} (\Pi_{+}) \cap C^{1} (\ol{\Pi}_{+})
\label{2.43}\\
\frac{\partial v}{\partial x_{n}} &= - f, \quad x \in
\partial \Pi_{+} \label{2.44}\\
v(x) &= o \left( | x |^{M} \sec^{n-2} \theta
\right); \quad x \in \Pi_{+}, \quad |x| \raro \infty .
\label{2.45}
\end{align}
\end{cor}

\bigskip
\noindent
{\it Proof}:  The growth estimate \eqref{2.45} follows
from the Theorem:
$$
N_{M} [wf] (x) = \frac{\alpha_{n}}{n-2}
F_{\frac{n-2}{2}, M} [wf](x).
$$
And,
\begin{align*}
\left| N[(1-w) f] (x) \right| & \leq
\frac{\alpha_{n}}{n-2} \int\limits_{| y' | < 2}
| f (y) |\, ( | x | -2 )^{2-n}\, d y' \\
& \leq \alpha_{n} (n-2)^{-1}\, 2^{n-2} | x |^{2-n}
\quad \text{if} \quad | x | \geq 4.
\end{align*}

Theorem 1 of \cite{Gardiner} shows \eqref{2.40},
\eqref{2.43} and \eqref{2.44} hold. \qed 

\begin{rem}
\label{rem2.3}{\rm  In Corollary \ref{cor2.1}, the
solution to \eqref{2.39}--\eqref{2.42} is unique if $M=0$
and if $M \geq 1$ it is unique to the addition of a
harmonic polynomial of degree $M$ vanishing on $\partial
\Pi_{+}$ (\cite{SiegelTalvila}, Theorem 3.1).
Similarly, in Corollary \ref{cor2.2}, if $M=0$ the
solution to \eqref{2.40}, \eqref{2.43}--\eqref{2.45} is unique and if
$M \geq 1$ it is unique to the addition of a harmonic
polynomial $p(x)$ of degree $M-1$ that is even about
$x_{n} = 0$.  This can be proved by modifying the proof
of Theorem 3.1 (\cite{SiegelTalvila}):  For a
function $v$ that satisfies \eqref{2.43}--\eqref{2.45}
with $f=0$, its even extension is defined in $\mathbb{R}^{n}$,
and by expanding $v$ in terms of spherical harmonics,  using
\eqref{2.45}, the conclusion is obtained.} \hfill
\end{rem}

\begin{rem}
\label{rem2.4}{\rm  If $f(y')\,|y'|^{n-2-M}$ is integrable
at the origin then we can use $u(x) = D_{M}[f](x)$ and
$v(x) = N_{M}[f](x)$ in Corollaries \ref{cor2.1} and
\ref{cor2.2}, respectively.}
\hfill

\end{rem}

\begin{cor}\label{cor2.3}
If $\omega\! : \Pi_{+} \raro (0, \infty )$ then $\omega$
is a sharp growth condition for $F_{\lambda , M}$ if and
only if there are constants $0 < S < T < \infty$ and
$N>0$ such that $S \leq | x |^{-M} \cos^{2
\lambda} \theta\, \omega (x) \leq T$ for all $x \in
\Pi_{+}$ with $| x | > N$.
\end{cor}

\bigskip
\noindent
{\it Proof}: Throughout the proof $f$ will satisfy
\eqref{1.5} and $| x |$, $| x^{(i)} | > N$.

Suppose $S$ and $T$ exist as above. Then
$$
\left| F_{\lambda , M} [f](x) \right| / \omega (x) \leq
\left| F_{\lambda , M} [f](x) \right| S^{-1} | x
|^{-M} \cos^{2 \lambda} \theta \raro 0
$$
so $F_{\lambda , M} [f] = o(\omega )$.

Let $\psi \!: \Pi_{+} \raro (0,\infty)$ with $\psi = o
(\omega )$ then $\psi (x) = o ( | x |^{M} \sec^{2
\lambda} \theta )$.  Given $\{ x^{(i)} \}$ in $\Pi_{+}$
take $f$ as in the proof of the Theorem.  Then
$F_{\lambda , M} [f](x^{(i)})/ \psi (x^{(i)}) \not\raro 0$.
Hence, $\omega$ is sharp.

Now suppose $\omega$ is sharp.  If $\chi (x) := |
x|^{-M} \cos^{2\lambda} \theta\, \omega (x)$ is
unbounded then there is a sequence $\{ x^{(i)} \}$ on
which $\chi (x^{(i)}) \raro \infty$.  Let
\begin{equation*}
\psi (x) = \begin{cases}
| x |^{M} \sec^{2\lambda} \theta & \text{on} \quad
\{ x^{(i)} \} \\
\omega (x) / | x |, & \text{otherwise,} \end{cases}
\end{equation*}
then $\psi = o(\omega )$ but for any  $f$ we have
$
F_{\lambda, M}[f](x^{(i)})/\psi (x^{(i)}) \raro 0
$,
since\newline
$
F_{\lambda ,M} [f](x) = o (| x |^{M}
\sec^{2\lambda} \theta ).
$
This contradicts the assumption that $\omega$ was sharp
(Definition 2.1, (ii)).  Hence $T$ exists as above.

If $\chi \raro 0$ on some sequence $\{ x^{(i)} \}$ then
take $f$ such that
\begin{equation}
\limsup_{i \raro \infty}\, F_{\lambda, M} [f](x^{(i)}) | x^{(i)}
|^{-M} \cos^{2\lambda} \theta_{i} \geq 1 \quad (\cos
\theta_{i} = x^{(i)}_{n}/| x^{(i)} |).
\end{equation}
Then
$$
\limsup_{i \raro \infty}\,
\frac{F_{\lambda , M} [f](x^{(i)})}{\omega (x^{(i)})} =
\limsup_{i \raro \infty}\,\frac{F_{\lambda , M}[f](x^{(i)})}{| x^{(i)} |^{M}
\sec^{2\lambda} \theta_{i}} \frac{| x^{(i)}|^{M}
\sec^{2\lambda} \theta_{i}}{\omega (x^{(i)})} =
\infty ,
$$
which contradicts the sharpness assumption (i) of
Definition 2.1.  Hence, $S$ exists as above. \qed 

\begin{rem}
\label{rem2.5}{\rm  The angular blow up predicted for
$F_{\lambda , M}$ as $| x | \raro \infty$ can be
expected to occur only as $x$ approaches $\partial
\Pi_{+}$ within a thin or rarefied set.  See
\cite{Aikawa}, \cite{Essen}, \cite{Mizuta}
and references therein.}\hfill
\end{rem}

\section{\!\!\!\!\!\!.\; Representation of the Neumann solution.}

The modified kernel $K_{M} (\lambda , x, y')$ satisfies a
differential-difference equation for the derivative with
respect to $\theta , | x | , y_{i}, | y | ,
x_{n}$, $| y' |$, $y_{i}'$ and $\theta'$, relating
the derivative to\newline $K_{M} (\lambda +1, x,y')$, $K_{M-1}
(\lambda +1, x, y')$ and $K_{M-2} (\lambda +1, x,y')$.
The integration of these equations give representations
of the modified Neumann integral in terms of the
modified Dirichlet integral.

\begin{prop}\label{prop3.1}
Let $n \geq 3$, $M\geq 0$, $\lambda > 0$, $x
\in \Pi_{+}$ and $y' \in \mathbb{R}^{n-1}$.  Use the convention
that $K_{m} =K$ if $m \leq 0$.  Then
\begin{itemize}
\item[(i)]
${\dss\frac{\partial K_{M}}{\partial \theta}} (\lambda ,
x, y') = 2 \lambda\, x_{n}\, \hat{y} \!\cdot\! y'\, K_{M-1}
(\lambda + 1, x, y')$

\item[(ii)]
${\dss\frac{\partial K_{M}}{\partial | x |}}
(\lambda , x, y') = 2 \lambda \left[ \sin \theta\,\, \hat{y}
\!\cdot\! y'\, K_{M-1} (\lambda +1, x, y') - | x |
K_{M-2} (\lambda +1, x, y') \right]$

\item[(iii)]
${\dss\frac{\partial K_{M}}{\partial y_{i}}} (\lambda ,
x, y') = 2 \lambda \left[ y_{i}'\, K_{M-1}\! (\lambda +1, x,
y') - y_{i} K_{M-2} (\lambda +1, x, y') \right],\!\!\!\! \quad\! {1
\leq i \leq n-1}$

\item[(iv)]
${\dss\frac{\partial K_{M}}{\partial | y |}}
(\lambda , x, y') = 2 \lambda \left[ \hat{y} \!\cdot\!
y'\, K_{M-1} (\lambda + 1,x,y') - | y |
K_{M-2} (\lambda +1 , x, y') \right]$

\item[(v)]
${\dss\frac{\partial K_{M}}{\partial x_{n}}} (\lambda ,
x, y') = -2 \lambda\, x_{n}\, K_{M-2} (\lambda +1, x, y')$

\item[(vi)]
${\dss\frac{\partial K_{M}}{\partial | y' |}}
(\lambda , x, y') = 2 \lambda \left[ y \!\cdot\! \hat{y}'\,
K_{M-1} (\lambda +1, x,y') - | y' | K_{M} (\lambda
+ 1, x, y') \right]$

\item[(vii)]
${\dss\frac{\partial K_{M}}{\partial y_{i}'}} (\lambda ,
x, y') = 2
\lambda \left[ y_{i}\, K_{M-1} (\lambda +1, x,y') - y_{i}'\,
K_{M} (\lambda +1,x,y') \right],\!\!\!\! \quad\! {1
\leq i \leq n-1}$

\item[(viii)]
${\dss\frac{\partial K_{M}}{\partial \theta'}} (\lambda ,
x, y') = - 2
\lambda\, | y | \, | y' | \sin \theta'\, K_{M-1}
(\lambda +1, x,y')$.
\end{itemize}
\end{prop}

The proofs rest on the identities
\begin{equation}
\frac{d}{dt} C_{m}^{\lambda} (t) = 2 \lambda\,
C_{m-1}^{\lambda +1} (t) \label{3.1}
\end{equation}
\begin{equation}
C_{0}^{\lambda} (t) =1, \quad 
C_{m}^{\lambda} (t) \equiv 0 \quad \text{for} \quad m =
0,-1, -2, \cdots  \label{3.2}
\end{equation}
\begin{equation}
mC_{m}^{\lambda}(t) = 2 \lambda [ t\, C_{m-1}^{\lambda +1}
(t) - C_{m-2}^{\lambda +1} (t)] \label{3.3}
\end{equation}
\begin{equation}
(m+2\lambda ) C_{m}^{\lambda} (t) = 2 \lambda
[C_{m}^{\lambda +1} (t) - t\, C_{m-1}^{\lambda +1} (t)]
\label{3.4}
\end{equation}
(\cite{Szego}, 4.7.28).  In (iv) $x_{n}$ is fixed.  If 
$\theta$ is held constant for the differentiation then
$${\dss\frac{\partial K_{M}}{\partial | y |}}
(\lambda , x, y') = 2 \lambda \left[ \hat{y} \!\cdot\!
y'\, K_{M-1} (\lambda + 1,x,y') - | x |\csc\theta\,
K_{M-2} (\lambda +1 , x, y') \right].$$  
This leads to a similar change in (iv)  of Proposition \ref{prop3.2}.

\bigskip
\noindent
{\it Proof of (i)}:  From \eqref{1.1}, \eqref{1.6},
\eqref{3.1} and \eqref{3.2}
\begin{align*}
\frac{\partial K_{M}}{\partial \theta} (\lambda , x, y')
&= 2 \lambda\, | y' | \, | x | \cos \theta
\cos \theta' K (\lambda +1, x, y') \\
& \qquad\qquad -2 \lambda \sum_{m=1}^{M-1} | x
|^{m} | y' |^{-(m+2 \lambda)} C_{m-1}^{\lambda +1}
(\Theta ) \cos \theta \cos \theta' \\
&= 2 \lambda\, x_{n}\, | y' | \cos \theta' K_{M-1}
(\lambda +1, x, y').\quad\,\, \blacksquare 
\end{align*}

The other proofs follow in a similar manner from
\eqref{3.1}--\eqref{3.4}.  Note that
$\hat{y}$ and $\Theta = \sin \theta\, \cos \theta' = \sin
\theta\,\, \hat{y} \!\cdot\! \hat{y}'$ are independent of $|
x|$ and that $\tan \theta = | y | / x_{n}$ so that
$\partial \theta/\partial y_{i} = y_{i}\,
x_{n}/(| x |^{2} | y | )$ and
$\partial \theta/\partial x_{n} = - \sin \theta
/ | x |$.

Now introduce the following notation.   If
$z_{1}$ and $z_{2}$ are in $\mathbb{R}^{n-1}$
then $f_{z_{1}} (z_{2}) = z_{1} \!\cdot\! z_{2}\, f(z_{2})$.  If
$0 \leq s \leq \pi /2$ then $x(s)$ indicates $x$
with the polar angle $\theta$ replaced by $s$, i.e.,
$x(s) = y(s) + x_{n} (s) \hat{e}_{n}$, where $y(s) =
| x | \sin s\, \hat{y}$, $x_{n} (s) = |x| \cos
s$, $x (\theta ) =x$ and $y(\theta ) = y$.  Note that
$| x |$ and $\hat{y}$ are independent of
$\theta$.  And, if $x =
{\dss\sum_{j=1}^{n}}x_{j} \hat{e}_{j}$ then
$\breve{x}_{i} (t) = {\dss\sum_{j\neq i}} x_{j}
\hat{e}_{j} + t \hat{e}_{i}$ $(1 \leq i \leq n)$.

Integrating (i) through (iv) above and setting $\lambda =
(n-2)/2$ we obtain

\begin{prop}\label{prop3.2} 
Let $f$ be continuous with the
origin not in the closure of its support and satisfy
\eqref{1.5} with $\lambda = (n-2)/2$ $(n \geq 3)$. Let
$M \geq 0$, $\lambda > 0$, $x \in \Pi_{+}$ and adopt the
convention that $D_{m} = D$ for $m \leq 0$.  Then the
following are equal to $N_{M} [f](x)$
\begin{itemize}
\item[(i)]
${\dss\int\limits_{t=\theta_{0}}^{\theta}} D_{M-1}
[f_{\hat{y}}] (x(t))\, dt + N_{M} [f] (x(\theta_{0})), \quad
0 \leq \theta_{0} \leq \frac{\pi}{2}$

\item[(ii)]
$\tan \theta {\dss\int\limits_{t=r_{0}}^{| x |}}
D_{M-1} [f_{\hat{y}}] (t \hat{y}) \frac{dt}{t} - \sec
\theta {\dss\int\limits_{t=r_{0}}^{| x |}} D_{M-2}
[f] (t \hat{x})\, dt 
+ N_{M} [f] (r_{0} \hat{x}), \quad r_{0} \geq 0$

\item[(iii)]
$\frac{1}{x_{n}} {\dss\int\limits_{t=t_{i}}^{y_{i}}}
D_{M-1} [f_{\hat{e}_{i}}] (\breve{x}_{i} (t))\, dt -
\frac{1}{x_{n}} {\dss\int\limits_{t=t_{i}}^{y_{i}}}
D_{M-2} [f] (\breve{x}_{i}(t))\, t\,dt$ \\
\begin{equation*}
\;\;\;\;\;\;\;\;\;\;\;\;\;\;\;\;\;+N_{M} [f] (\breve{x}_{i} (t_{i})); 
\quad t_{i} \in
\mathbb{R},\,1 \leq i \leq n -1
\end{equation*}

\item[(iv)]
$\frac{1}{x_{n}}{\dss\int\limits_{t=\rho}^{| y|}}
D_{M-1}[f_{\hat{y}}](t \hat{y} + x_{n}
\hat{e}_{n})\,dt -
\frac{1}{x_{n}}{\dss\int\limits_{t=\rho}^{| y |}} D_{M-2}
[f](t\hat{y} + x_{n} \hat{e}_{n})\,t\,dt
$ \\
\begin{equation*}
 + N_{M} [f] (\rho \hat{y} + x_{n} 
\hat{e}_{n}),\quad
\rho \geq 0
\end{equation*}

\item[(v)]
$- {\dss\int\limits_{t=t_{n}}^{x_{n}}} D_{M-2} [f]
(\breve{x}_{n} (t))\, dt + N_{M} [f] (\breve{x}_{n}
(t_{n})),\quad t_{n} \geq 0$.
\end{itemize}
\end{prop}

\bigskip
\noindent
{\it Proof}:  Integrate each of (i)--(v) in Proposition
\ref{prop3.1} with respect to the relevant variable and set
$\lambda = (n-2)/2$.  Multiply by $((n-2)/2)f(y')$ and
integrate $y' \in \mathbb{R}^{n-1}$.  Because of \eqref{1.5}
the integrals $D_{M-2} [| f | ] (x)$ and $D_{M-1}
[ | f_{\hat{y}} | ] (x)$ converge to continuous
functions on $\ol{\Pi}_{+}$.  The same is true for each
modified Dirichlet integral in (i)--(v).  Fubini's
Theorem now justifies the interchange of orders of
integration. \qed 

\begin{rem}
\label{rem3.1}{\rm We can relax the condition that $f$ be
continuous if we refrain from evaluating $N_{M}[f]$ on
$\partial \Pi_{+}$.  This requires taking $0 \leq
\theta_{0} < \pi /2$, $r_{0} >0$, $\rho > 0$ and $t_{n}
>0$.  We can dispense with the restriction on the support
of $f$ if $f(y')\,|y'|^{-(M+n-3)}$
is integrable at the origin.
When $M=0$ there is no restriction on the support of $f$.}
\end{rem}

We can use Proposition \ref{prop3.2}(i) to confirm the growth
estimate \eqref{2.45}.  If $f$ satisfies \eqref{1.5}
with $\lambda = (n-2)/2$, $n \geq 3$, and if $0 \leq
\theta_{0} < \pi /2$ then \eqref{4.2.17} gives $N_{M} [f]
(x(\theta_{0})) = o (| x |^{M})$.   From Corollary
\ref{cor2.1}, $D_{M-1} [f_{\hat{y}}] (x) = o (| x
|^{M} \sec^{n-1} \theta )$.  Integrating over
$\theta$ and using (i) of Proposition \ref{prop3.2},
$$
N_{M} [f](x) = o ( | x |^{M} \sec^{n-2} \theta ) +
o (| x |^{M}) = o (| x |^{M} \sec^{n-2}
\theta ) ,
$$
in agreement with Corollary \ref{cor2.2}.

\section{\!\!\!\!\!\!.\; Second type of modified kernel.}

Using the generating function \eqref{2.1} with $z = |
y ' | / | x |$, $t = \Theta = \sin \theta \cos
\theta'$, we can define a second type of modified kernel
\begin{equation}
\tilde{K}_{M} (\lambda , x, y') = K(\lambda , x, y') -
\sum_{m=0}^{M-1} \frac{| y' |^{m}}{| x
|^{m+2\lambda}} C_{m}^{\lambda} (\Theta ) ,
\label{4.1}
\end{equation}
defined for $| x | > 0$ and $M \geq 1$.  The
convergence condition corresponding to \eqref{1.5} is now
\begin{equation}
\int\limits_{\mathbb{R}^{n-1}} | f (y') | ( | y'
|^{M-1} +1)\, dy' < \infty . \label{4.2}
\end{equation}
If \eqref{4.2} is satisfied, define
\begin{equation}
\tilde{F}_{\lambda , M} [f] (x) = \int\limits_{\mathbb{R}^{n-1}}
f(y') \tilde{K}_{M} (\lambda , x, y')\, dy' . \label{4.3}
\end{equation}
Define $\tilde{D}_{M}$ and $\tilde{N}_{M}$ in terms of 
$\tilde{F}_{\lambda , M}$ as in
\eqref{1.7} and \eqref{1.8}.  Each $x_{n} |
x |^{-(m+n)} C_{m}^{n/2} (\Theta )$ in the
kernel $\tilde{D}_{M}$ is harmonic in $\mathbb{R}^{n}\backslash
\{ 0 \}$ (Remark \ref{rem2.2}).  The same can be said for $| x
|^{-(m+n-2)} C_{m}^{(n-2/2)} (\Theta )$ in the kernel
$\tilde{N}_{M}$.  Hence, $\tilde{D}_{M} [f]$ and
$\tilde{N}_{M} [f]$ are harmonic in $\Pi_{+}$.  Results
similar to Propositions \ref{prop3.1} and \ref{prop3.2} hold for
$\tilde{K}_{M}$, $\tilde{D}_{M}$ and $\tilde{N}_{M}$.

However, $\tilde{D}_{M} [f]$ is not continuous on
$\ol{\Pi}_{+}$.  Since \eqref{4.2} implies \eqref{1.4}
(with $\lambda = n/2$) the unmodified Poisson integral
$D[f]$ is continuous for $x_{n} \geq 0$ if $f$ is
continuous.  Hence,
$$
\tilde{D}_{M} [f] (x) = D[f](x) - \alpha_{n} x_{n}
\sum_{m=0}^{M-1} | x |^{-(m+n)}
\int\limits_{\mathbb{R}^{n-1}} | y' |^{m} f(y')
C_{m}^{n/2} (\Theta )\, dy'
$$
and $\tilde{D}_{M}[f]$ is continuous for $x_{n} \geq 0$,
$x \neq 0$.  Similar remarks apply to the Neumann case.
We will work with $\tilde{D}_{M}$ and $\tilde{N}_{M}$
only in the limit $| x | \raro \infty$.

Growth estimates for $\tilde{F}_{\lambda , M}$ are
similar to those for $F_{\lambda , M}$.

\begin{theo}\label{theo4.1}
If \eqref{4.2} holds for measurable $f$ then
$$
\tilde{F}_{\lambda , M} [f] (x) = o (| x
|^{-(M+2\lambda -1)} \sec^{2\lambda} \theta )\quad (x \in
\Pi_{+},\,| x | \raro \infty )
$$
and this estimate is sharp in the sense of Definition 
\ref{defin2.1}.
\end{theo}

\bigskip
\noindent

{\it Proof}: Throughout the proof $d_{1}$ and $ d_{2}$
will be positive constants (depending on
$\lambda$ and $M$).  In \eqref{2.11} replace $| x
|/ | y' |$ by $| y' | / | x |$ and
in the proof of Theorem \ref{theo2.1} let $s = | y'
| / | x |$.

If $0 < \lambda < \frac{1}{2}$ then \eqref{4.2.14},
\eqref{4.2.15} and \eqref{4.2.18} give
\begin{align*}
| \tilde{K}_{M} (\lambda , x, y') | & \leq d_{1}\,
K(\lambda , x, y')\, s^{M-1} \int\limits_{\zeta =0}^{s}
| 1- \zeta |^{2\lambda -1}\, d \zeta \\
& \leq \frac{d_{2}\, s^{M+2 \lambda -1}
\sec^{2\lambda}\theta}{(| x | + | y' |
)^{2\lambda}} ,
\end{align*}
from which $\tilde{F}_{\lambda , M}[f](x) = o (| x
|^{-(M+2 \lambda -1)} \sec^{2\lambda} \theta )$.

If $\lambda \geq \frac{1}{2}$ then
\eqref{4.2.14} and \eqref{4.2.15} give
\begin{align}
| \tilde{K}_{M} (\lambda , x, y')| &
\leq 
\f{d_2|y'|^{M-1}\sec^{2\la}\theta}{| x |^{M+2\la -1}}\left[
1-(1+s)^{-2\lambda}\right]
\label{5.4.4}\\
 & \leq \f{d_2|y'|^{M-1}\sec^{2\la}\theta}{| x |^{M+2\la -1}}.
\notag
\end{align}
Since $1-(1+s)^{-2\la}\to0$ as $s\to 0$, integrating \eqref{5.4.4}
and noting \eqref{4.2}, dominated convergence gives
$$
\int\limits_{\Rbn} f(y')
\tilde{K}_{M} (\lambda , x, y')\, dy' = o (| x
|^{-(M+2\lambda -1)} \sec^{2\lambda} \theta )\quad (x \in
\Pi_{+},\, | x | \raro \infty ) .
$$

To prove this sharp, interchange $| x |$ and $|
y' |$ in the proof of Theorem \ref{theo2.1} and
proceed in a similar manner.  \qed

The modified kernel furnishes an asymptotic expansion of
$D[f]$ and $N[f]$.

\begin{theo}\label{theo4.2}
Let $f$ be measurable such that \eqref{4.2} holds for a
positive integer $M$.  Then, as $x\raro \infty$ in
$\Pi_{+}$
\begin{itemize}
\item[(i)]
$D[f](x) = {\dss\sum_{m=0}^{M-1}} | x |^{-(m+n-1)}
Y_{m+1}^{(0)} (\hat{x}) + o (| x |^{-(M+n-2)}
\sec^{n-1}\theta)\quad (n \geq 2)$

\item[(ii)]
$N[f](x) = {\dss\sum_{m=0}^{M-1}} | x |^{-(m+n-2)}
Y_{m}^{(1)} (\hat{x}) + o (| x |^{-(M+n-3)}
\sec^{n-2} \theta )\quad (n \geq 3)$
\end{itemize}
where $Y_{m}^{(0)}$ is given by \eqref{4.5} below and is
a spherical harmonic of degree $m$ that vanishes on
$\partial \Pi_{+}$ and $Y_{m}^{(1)}$ is given by
\eqref{4.6} below and is a spherical harmonic of degree
$m$ whose normal derivative vanishes on $\partial
\Pi_{+}$.  The data $f$ can be chosen so that
simultaneously the leading order term $(m=0)$ does not
vanish and the order relation is sharp (in the sense of
Definition \ref{defin2.1}).
\end{theo}

\bigskip
\noindent
{\it Proof}:  To prove (i) use \eqref{4.1} with $\lambda
= n/2$, $f$ as in the Theorem and $| x | > 0$,
$$
D[f](x) = \alpha_{n} x_{n} \int\limits_{\mathbb{R}^{n-1}} f(y')
\sum_{m=0}^{M-1} \frac{| y' |^{m}}{| x
|^{m+n}} C_{m}^{n/2} (\Theta )\,dy' + \tilde{D}_{M} [f](x).
$$
Now (i) follows from Theorem \ref{theo4.1} and the
definition
\begin{equation}
Y_{m+1}^{(0)} (\hat{x}) = \alpha_{n} \cos \theta
\int\limits_{\mathbb{R}^{n-1}} f(y')\, | y' |^{m}\, C_{m}^{n/2}
(\sin \theta\,\, \hat{y} \!\cdot\! \hat{y}')\, dy' . \label{4.5}
\end{equation}

Clearly $Y_{m+1}^{(0)}$ vanishes when $\theta = \pi
/2$.  It is a spherical harmonic of degree $m+1$ by
Remark \ref{rem2.2}.

Given $\psi (x) = o ( | x |^{-(M+n-2)}\sec^{n-1}
\theta )$ take $f$ as in Theorem \ref{theo4.1} so that
$\tilde{D}_{M} [f] (x) = o (| x |^{-(M+n-2)}
\sec^{n-1} \theta )$ is sharp.  In particular, $f$ can
be taken to be positive for $y_{1} > 0$ with a super-even
extension ($M$ even) or super-odd extension ($M$ odd)
to $y_{1} < 0$ ($A_{\lambda} > 1$ in the proof of
Theorem \ref{theo2.1} \eqref{4.2.32}.   See Step III in the outline
of the proof for an explanation of the terminology ``super-even" 
and ``super-odd.")
The leading order
term in (i) is with $m=0$, $Y_{1}^{(0)} (\hat{x}) =
\alpha_{n} \cos \theta \int\limits_{\mathbb{R}^{n-1}} f(y')\,
dy'$.  With $f$ as above this spherical harmonic does
not vanish if $0 \leq \theta < \pi /2$.
With the definition
\begin{equation}
Y_{m}^{(1)} (\hat{x}) = \frac{\alpha_{n}}{n-2}
\int\limits_{\mathbb{R}^{n-1}} f(y')\, | y' |^{m}\,
C_{m}^{(n-2)/2} (\Theta )\, dy',\label{4.6}
\end{equation}
the proof of (ii) is similar.\qed 

The addition formula for Gegenbauer polynomials can be
used to separate the $\theta$ dependence in (i) and
(ii).  First write
$$
C_{m}^{n/2} (\sin \theta\,\, \hat{y} \!\cdot\! \hat{y}') =
\sum_{\ell =0}^{\left\lfloor \frac{m}{2} \right\rfloor}
\gamma_{n,m,\ell} (\theta\, ) C_{m-2\ell}^{(n-1)/2}
(\hat{y} \!\cdot\! \hat{y}'),
$$
where
\begin{equation*}
\begin{split}
\gamma_{n,m,\ell} (\theta ) &= \frac{(n-2)!\, (-1)^{\ell}
\, (2\ell )!\, (n+2m-4\ell -1)\, \Gamma
\left(\frac{n}{2} +m-2\ell \right)\, \Gamma \left(
\frac{n}{2} +m-\ell \right)}{4^{2\ell -m}\,\Gamma^{2}
(n/2)\, \ell !\, (n+2 m-2\ell -1)!}\times \\
& \;\;\;\;\;\;\;\;\;\;\;\;\;\;\qquad \qquad \times \sin^{m-2\ell} \theta\,
C_{2\ell}^{n/2 +m-2\ell} (\cos \theta )
\end{split}
\end{equation*}
\renewcommand{\thefootnote}{\fnsymbol{footnote}}
(\cite{Erdelyi} 10.9.34\footnote[2]{The first term in the sum
over $m$ in this formula should read $2^{2m}$.},
10.9.19).  Then \eqref{4.5} becomes
$$
Y_{m+1}^{(0)} (\hat{x}) = \alpha_{n} \cos \theta
\sum_{\ell = 0}^{\left \lfloor \frac{m}{2}
\right\rfloor} \gamma_{n,m,\ell} (\theta )\, \delta_{n, m,
\ell} (\hat{y}),
$$
where
$$
\delta_{n,m,\ell} (\hat{y}) = \int\limits_{\mathbb{R}^{n-1}}
f(y')\, | y' |^{m}\, C_{m-2\ell}^{(n-1)/2} (\hat{y}
\!\cdot\!\hat{y}')\, dy'
$$
and is independent of $| x |$ and $\theta$.

A similar separation of $| x |$, $\theta$ and
$\hat{y}$ dependence in (ii) is given by
$$
Y_{m}^{(1)} (\hat{x}) = \frac{\alpha_{n}}{n-2}
\sum_{\ell =0}^{\left\lfloor \frac{m}{2} \right\rfloor}
\gamma_{n-2, m, \ell} (\theta )\, \delta_{n-2, m, \ell}
(\hat{y}) \quad (n \geq 3).
$$

The function defined by $Z_{m} (\hat{y}_{1} ,
\hat{y}_{2}) = C_{m}^{n/2} (\hat{y}_{1} \!\cdot\!
\hat{y}_{2})$ is known as a zonal harmonic of degree $m$
with pole $\hat{y}_{1} \in \partial B_{1}$, evaluated at
$\hat{y}_{2} \in \partial B_{1}$ (see \cite{Axler},
Chapter 5).

If the integral in \eqref{4.2} converges for all $M
\geq 1$, letting $M \raro \infty$ in Theorem
\ref{theo4.2} will give asymptotic series for $D_{M}[f]$
and $N_{M} [f]$.  As the following example shows, these
series will not in general be convergent.

\bigskip
\noindent{\bf Example.}
Let $f(y) = \exp (- | y |)$ and let $d \omega_{n-1}$ be
surface measure on the unit ball of $\mathbb{R}^{n-1}$.  Then for
$n \geq 3$, \eqref{4.6} becomes
\begin{align*}
Y_{m}^{(1)} (\hat{x}) &= \frac{\alpha_{n}}{n-2}
\int\limits_{\rho =0}^{\infty} e^{-\rho} \rho^{m+n-2}\, d
\rho \int\limits_{\partial B_{1}} C_{m}^{(n-2)/2} (\sin
\theta\,\, \hat{y} \!\cdot\! \hat{y}')\, d\omega_{n-1} \\
& = \frac{\alpha_{n}}{n-2}\, (m+n-2)!\, (n-2)\, \omega_{n-2}\,
I_{n,m}^{(1)} (\theta ),
\end{align*}
where
$$
I_{n,m}^{(1)} (\theta ) = \int\limits_{\phi=0}^{\pi}
C_{m}^{(n-2)/2} (\sin \theta \cos \phi) \sin^{n-3}\phi\, d \phi
$$
and the surface integral was evaluated by spherical
means \cite{John}.  The integral $I_{n,m}^{(1)} (\theta
) $ is known (\cite{GradshteynReyzhik} 7.323.2, together
with \cite{Erdelyi} 10.9.19),
\begin{equation*}
I_{n,m}^{(1)} (\theta ) = \begin{cases}
\frac{2^{n-3} \Gamma( n/2\, -1)\,
(-1)^{k}\, (2k)!\, \Gamma ( k + n/2\, -1 )\,
C_{2k}^{(n-2)/2} (\cos \theta )}{k!\, \Gamma (2k
+n-2)},\quad m=2k \\
0, \quad m \quad \text{odd}. \end{cases}
\end{equation*}
And,
$$
Y_{2k}^{(1)} (\hat{x}) = \frac{2^{n-2} \Gamma( n/2 -1)\,
(-1)^{k}\, (2k)!\, \Gamma ( k + n/2 )\,
C_{2k}^{(n-2)/2} (\cos \theta )}
{\pi k!}
$$
$(Y_{2k+1}^{(1)} (\hat{x}) =0)$.  As $k \raro \infty$,
$$
C_{2k}^{(n-2)/2} (\cos \theta ) \sim \frac{2\left(
n/2 +2k -2 \right) !\, \sin \left[ n\pi/4
- \left( n/2 +2k-1 \right) \theta
\right]}{\left( n/2 -2 \right) !\, (2k)!\, (2 \sin
\theta )^{n/2-1}}
$$
(\cite{Szego} 8.4.13), so Stirling's approximation shows that for
fixed $x$
\begin{equation*}
\sum_{m=0}^{M-1} | x |^{-m} Y_{m}^{(1)} (\hat{x})
\quad \text{diverges as} \quad M \raro \infty .
\end{equation*}
With the Dirichlet expansion we have from \eqref{4.5}
$$
Y_{m+1}^{(0)} (\hat{x}) = \alpha_{n} (m+n-2)!\, (n-2)\,
\omega_{n-2}\, I_{n,m}^{(0)} (\theta ) ,
$$
where
$$
I_{n,m}^{(0)} (\theta ) = \cos \theta
\int\limits_{\phi=0}^{\pi} C_{m}^{n/2} (\sin \theta \cos \phi)
\sin^{n-3} \phi\,  d\phi.
$$
If $n \geq 5$ then \eqref{2.3} and integration by parts
give
$$
I_{n,m}^{(0)} (\theta ) =\frac{1}{(n-2) \sin \theta}\,
\frac{d}{d\theta} I_{n-2, m+2}^{(1)} (\theta )
$$
and
\begin{equation}
\sum_{m=0}^{M-1} | x |^{-m} Y_{m}^{(1)} (\hat{x})
\quad \text{diverges as} \quad M \raro \infty .\label{4.7}
\end{equation}

When $n=2$, we use
$$
C_{m}^{1} (\cos \phi) = \frac{\sin [(m+1)\phi]}{\sin \phi}
$$
and replace $\theta$ by $\pi /2 -\phi$ (see the end of the proof
of Lemma  \ref{lem2.1}).
With Dirichlet data $f(\xi ) = \exp (- | \xi |
)$ we have $Y_{m+1}^{(0)} (\hat{x}) = 2\,m! \sin
[(m+1)\phi]/\pi$,  if $m$ is even  and again
$$
\sum_{m=0}^{M-1} \frac{m! \sin [(m+1)\theta ]}{r^{m}}
\quad \text{diverges as} \quad M \raro \infty
$$
for fixed $x = r\cos \theta\,\hat{e}_{1} + r\sin\theta\,\hat{e}_{2} 
$.

When $n=3$, $I_{3,m}^{(0)}(\theta)$ can be evaluated in terms of Legendre
polynomials (since $C_{m}^{1/2}(t) = P_{m}(t)$) and when $n=4$,
$I_{4,m}^{(0)}(\theta)$ can be evaluated in terms of trigonometric functions.
In both  cases, the conclusion of \eqref{4.7} remains valid.

\end{document}